\documentclass[reqno,a4paper,final]{amsart}
\usepackage{amssymb,amstext,amsthm,eucal}
\usepackage{bbold}
\usepackage{array}
\usepackage[latin1]{inputenc}
\newcommand{\ecke}{\;_-\!\rule{0.2mm}{0.2cm}\;\;}

\newcommand{\SO}{\mathrm{SO}}

\newcommand{\Spin}{\mathrm{Spin}}

\newcommand{\GL}{\mathrm{GL}}

\newcommand{\SU}{\mathrm{SU}}
\newcommand{\U}{\mathrm{U}}

\newcommand{\RR}{\mathbb{R}}
\newcommand{\CC}{\mathbb{C}}

\newcommand{\1}{\mathbb{1}}

\newtheorem{THEO}{\bf Theorem}

\newtheorem{PR}{\bf Proposition}

\newtheorem{CO}{\bf Corollary}
\newtheorem{LM}{\bf Lemma}

\newcommand{\MUNCH}[1]{\relax}
\allowdisplaybreaks[1]
\begin{document}
\begin{sloppypar}
\title{Normal conformal Killing forms}
\author{Felipe Leitner}
\address{Institut f{\"u}r Mathematik, Universit{\"a}t Leipzig, Germany}
\email{leitner@mathematik.uni-leipzig.de}
\thanks{ }
\date{\today}

\begin{abstract} We introduce in this paper normal twistor equations for
differential forms and study their solutions, the
so-called normal conformal Killing forms.  The twistor equations arise
naturally from the canonical normal Cartan connection of conformal
geometry. Reductions of its holonomy are related to
solutions of the normal twistor equations. The case of decomposable normal
conformal holonomy representations is discussed.
A typical example with
an irreducible holonomy representation are the so-called
Fefferman spaces. We also apply our results to describe the geometry
of
solutions with conformal Killing spinors on Lorentzian spin
manifolds.

\end{abstract}
\maketitle
\tableofcontents

%%%%%%%%%%%%%%%%%%%%%%%%%%%%%%%%
\section{Introduction}%%%
%%%%%%%%%%%%%%%%%%%%%%%%%%%%%%%%
\label{p1}

A classical object of interest in differential geometry
is conformal symmetry.
Typical examples for conformal symmetry arise from the flow
of Killing
and conformal Killing vector fields on a 
semi-Riemannian manifold. The notion of conformal vector fields has
a natural generalisation to differential
forms and spinor fields, namely the
so-called conformal Killing forms and spinors
(cf. \cite{Kas68}, \cite{Tac69},\cite{Pen67}, \cite{PR86}, \cite{Lic88}, 
\cite{BFGK91},
\cite{Sem01}).

We want to introduce in this paper a special class
of conformal Killing forms, which we
call the normal conformal Killing forms (shortly:
nc-Killing forms).
These objects are solutions of certain
twistor equations, which are conformally covariant,
and moreover, they are subject to a normalisation
condition.
Their existence reflects a special part of the
conformal symmetry for
a metric (or conformal structure) on a semi-Riemannian
manifold.

The normal twistor equations are induced
by the canonical normal Cartan connection 
of conformal geometry.
This canonical connection is 
a well-known object in conformal geometry (cf. \cite{Kob72}). 
It lives on the principal
fibre bundle, which has a parabolic
subgroup of the M{\"o}bius group
as structure group.
Thereby, the parabolic subgroup consists of
those conformal transformations of the 
conformally flat model (M{\"o}bius space),
which fix the point at infinity. 
The normalisation condition for the normal
conformal Cartan connection is expressed in form of the 
$\varrho$-tensor
\[K_g=\frac{1}{n-2}\big(\frac{scal_g}{2(n-1)}-Ric_g\big)\ ,\]
which is a basic curvature tensor in conformal geometry (here given in 
terms of a metric 
$g$). 

The normal conformal Cartan connection can be
naturally extended to a usual principal fibre bundle
connection
with the M{\"o}bius group as structure group. 
Using this extended normal conformal connection, 
the notion of normal holonomy of a 
conformal structure can be introduced.
It is the structure group to which the normal conformal
connection can be `maximally' reduced. Since
nc-Killing forms find its interpretation as
parallel sections in certain tractor bundles
with respect to the covariant
derivative induced by the extended normal conformal
connection, the existence of nc-Killing forms
is apparently related to the holonomy group of the normal 
conformal connection.  

Similar as for the holonomy theory of the Levi-Civita
connection in (semi)-Riemannian geometry, the holonomy
of the normal conformal connection can be used    
for the characterisation of the underlying conformal geometry.
It turns out in the course of our discussion that a decomposable normal 
conformal holonomy
representation is related to the conformal Einstein
condition on (semi)-Riemannian spaces furnished with classical geometric 
structures such as Sasaki structures, nearly-K{\"a}hler structures 
etc.
In particular, it is 
possible to relate decomposability
of the normal conformal holonomy representation to the 
existence  
of a certain product metric
in the conformal class of a space.
The irreducible holonomy representations forbid the conformal Einstein
condition. A well-known example where this happens are the 
Fefferman spaces
in Lorentzian geometry, which arise by construction from CR-geometries.
The approach via the holonomy discussion of the normal conformal 
connection
will allow us to derive a certain geometric description for conformal
spaces admitting solutions of the normal twistor equations.
In particular, we will be able to give an improved geometric
characterisation of conformal spin spaces admitting 
conformal Killing spinors (cf. \cite{BL03}). 

The road map for our investigations of nc-Killing forms
and normal conformal holonomy is as follows.
In the paragraphs \ref{p2} to \ref{p4} we develop 
the basic
notions and facts for the construction of the canonical normal 
connection of conformal geometry and present finally the normal
twistor equations for differential forms.
In paragraph
\ref{p5} we derive integrability conditions
in terms of curvature conditions for the existence of
solutions. In paragraph
\ref{p6} we study solutions for nc-Killing forms on Einstein spaces
(cf. Theorem \ref{TH1}).
As we will
see, this is a natural thing to do in view
of the normality condition in form of the 
$\varrho$-tensor $K$.
In paragraph \ref{p7} we discuss the simplest
form of solutions, the
decomposable twistors, and we will understand in  
paragraph \ref{p8} that the 
solutions
on Einstein spaces are the basic building blocks, 
which appear
for the decomposable case (Theorem \ref{TH2}).
In paragraph \ref{p9}, we discuss as a showcase the normal
conformal holonomy representations and solutions of the twistor equations
on $4$-dimensional Riemannian and Lorentzian manifolds.
Finally, we use our results to discuss the conformal
Killing spinor equation on Lorentzian spin manifolds (cf. 
paragraph \ref{p10}; Theorem \ref{TH5}).

%%%%%%%%%%%%%%%%%%%%%%%%%%%%%%%%
\section{The representations}%%%
%%%%%%%%%%%%%%%%%%%%%%%%%%%%%%%%
\label{p2}
Let $\RR^{r,s}$ denote the (pseudo)-Euclidean space 
of signature $(r,s)$ with dimension $n=r+s$. 
The Lie algebra of conformal
Killing vector fields on $\RR^{r,s}$ is isomorphic to
$\frak{so}(r+1,s+1)$.
We explain here the usual action of 
$\frak{so}(r+1,s+1)$ on the spaces of $p$-forms 
$\Lambda^{p}_{r+1,s+1}$ over $\RR^{r+1,s+1}$
in terms of $2$-forms with respect to the irreducible
parts of the subrepresentations belonging to
the subalgebra $\frak{so}(r,s)$. The latter one
is the Lie algebra of the special orthogonal 
group $\SO(r,s)$,  
which is isomorphic to the set of 
Killing vector fields on
$\RR^{r,s}$ having a zero at the origin, i.e.,  
these are generators of orthogonal rotations.
In the following, we denote by $\flat$ and 
$\sharp$ the mappings between $\RR^{r,s}$ and its dual $\RR^{r,s*}$, which
are
naturally induced via the metric product $\langle\cdot,\cdot\rangle$ on
$\RR^{r,s}$. Moreover, we denote by $e=(e_1,\ldots,e_n)$ 
the standard orthonormal basis in $\RR^{r,s}$ such that  
$\varepsilon_i:=\langle e_i,e_i\rangle=-1$ for $i<r+1$.

The space of $2$-forms on $\RR^{r,s}$ is naturally isomorphic to
$\frak{so}(r,s)$ via the mapping
\[\begin{array}{cccc} 
\iota:&\Lambda^2_{r,s}&\rightarrow&\frak{so}(r,s)\subset\frak{gl}(n)\ 
.\\[2mm]
&\omega&\mapsto& (\ x\mapsto(x\ecke\omega)^\sharp\ )\end{array}\]
The natural action of $\Lambda^2_{r,s}$ on $\alpha\in\Lambda^p_{r,s}$ is 
then given by 
\begin{eqnarray*} e_i^\flat\wedge e_j^\flat\circ\alpha&=&
-e_i^\flat\wedge(e_j\ecke\alpha)+e_j^\flat\wedge(e_i\ecke\alpha)\\
&=&\ \  e_j\ecke(e_i^\flat\wedge\alpha)-e_i\ecke(e_j^\flat\wedge\alpha)\ .
\end{eqnarray*}

The Lie algebra $\frak{so}(r+1,s+1)$ of the group of
conformal transformations on the conformal compactification space 
$S^{r,s}$
of $\RR^{r,s}$ (M{\"o}bius space of signature $(r,s)$) is $|1|$-graded: 
\[\frak{so}(r+1,s+1)=\frak{g}_-\oplus\frak{g}_0\oplus\frak{g}_+\ ,\]
where $\frak{g}_-\cong\RR^{r,s}$, 
$\frak{g}_0\cong\frak{co}(r,s)$
and $\frak{g}_+\cong\RR^{r,s*}$ (see their brackets below). To set up
explicit 
identifications for these three subspaces,
let $(e_t,e_s,e_1,\ldots,e_n)$ be an orthonormal frame
of $\RR^{r+1,s+1}$, where $e_t$ is timelike, $e_s$ spacelike and
the $e_i$'s are the basis of $\RR^{r,s}$. We denote 
$e_-=\frac{1}{\sqrt{2}}(e_s-e_t)$ and 
$e_+=\frac{1}{\sqrt{2}}(e_s+e_t)$. Then we identify
\[\begin{array}{cccc}\iota:&\RR^{r,s}
&\rightarrow&\frak{g}_-\ , \\
&x&\mapsto& e_-^\flat\wedge x^\flat\end{array}\qquad\ \
\begin{array}{cccc}\iota:&\RR^{r,s*}&\rightarrow&\frak{g}_+\ , \\
&y^\flat&\mapsto& e_+^\flat\wedge y^\flat\end{array}\]
\[\begin{array}{cccc}\iota:&\RR\oplus\frak{so}(r,s)&
\rightarrow&\frak{g}_0\ .
\\ &(l,\omega)&\mapsto& l\cdot e_-^\flat\wedge
e_+^\flat+\omega\end{array}\]
Besides the usual bracket on
$\frak{g}_0\cong\frak{co}(r,s)$, the non-vanishing Lie
brackets are 
\[ [\omega,x]=(x\ecke\omega)^\sharp,\qquad [\omega,y^\flat]=
y\ecke\omega\quad\ \mbox{and}\quad\
[x,y^\flat]=\langle x,y\rangle\cdot e_-^\flat\wedge 
e_+^\flat+x^\flat\wedge
y^\flat\ ,\] 
where $x\in\frak{g}_-$, $y^\flat\in\frak{g}_+$ and
$\omega\in\frak{g}_0$. The
brackets
$[\frak{g}_-,\frak{g}_-]$ and $[\frak{g}_+,\frak{g}_+]$ all
vanish.

An arbitrary $(p+1)$-form $\alpha\in\Lambda^{p+1}_{r+1,s+1}$ 
on $\RR^{r+1,s+1}$ decomposes 
into
\[\alpha=e_-^\flat\wedge \alpha_-+\alpha_0+e_-^\flat\wedge
e_+^\flat\wedge\alpha_{\mp}+e_+^\flat\wedge\alpha_+\]
with uniquely determined forms $\alpha_-,\alpha_+\in\Lambda^{p}_{r,s}$,
$\alpha_0\in\Lambda^{p+1}_{r,s}$
and $\alpha_{\mp}\in\Lambda^{p-1}_{r,s}$. This split 
sum is with respect to the 
decomposition of $\Lambda^{p+1}_{r+1,s+1}$
into the irreducible submodules 
\[\Lambda^p_{r,s}\oplus\Lambda^{p+1}_{r,s}\oplus\Lambda^{p-1}_{r,s}\oplus
\Lambda^{p}_{r,s}\]
of the restricted action to $\frak{so}(r,s)$. 
The action of
$\frak{so}(r+1,s+1)$ 
on 
$\Lambda^{p+1}_{r+1,s+1}$ with respect to this decomposition
is given by
\[\begin{array}{lcl}e_-^\flat\wedge e_i^\flat\circ\alpha_-&=&0\\[1mm]
e_-^\flat\wedge
e_i^\flat\circ\alpha_0&=&-e_-^\flat\wedge(e_i\ecke\alpha_0)\\[1mm]
e_-^\flat\wedge e_i^\flat\circ\alpha_{\mp}&=&e_-^\flat\wedge
e_i^\flat\wedge
\alpha_{\mp}\\[1mm]
e_-^\flat\wedge
e_i^\flat\circ\alpha_+&=&e_i^\flat\wedge\alpha_++e_-^\flat\wedge e_+^\flat
\wedge(e_i\ecke\alpha_+)\end{array}\]
for $e_-^\flat\wedge e_i^\flat\in\frak{g}_-$. For 
$e_+^\flat\wedge e_i^\flat\in\frak{g}_+$ we have 
\[\begin{array}{lcl}
e_+^\flat\wedge
e_i^\flat\circ\alpha_-&=&e_i^\flat\wedge\alpha_--e_-^\flat\wedge e_+^\flat
\wedge(e_i\ecke\alpha_-)\\[1mm]e_+^\flat\wedge
e_i^\flat\circ\alpha_0&=&-e_+^\flat\wedge(e_i\ecke\alpha_0)\\[1mm]
e_+^\flat\wedge e_i^\flat\circ\alpha_{\mp}&=&-e_+^\flat\wedge
e_i^\flat\wedge
\alpha_{\mp}\\[1mm]
e_+^\flat\wedge e_i^\flat\circ\alpha_+&=&0\ \end{array}\]
and it is \[e_-^\flat\wedge e_+^\flat\circ\alpha=-e_-^\flat\wedge
\alpha_-+e_+^\flat\wedge\alpha_+\ .\]  
The action of $\frak{co}(r,s)$ on the components of $\alpha$ is the usual
one.

%%%%%%%%%%%%%%%%%%%%%%%%%%%%%%%%%%%%%%%%%%%%
\section{The normal conformal connection}%%%
%%%%%%%%%%%%%%%%%%%%%%%%%%%%%%%%%%%%%%%%%%%%
\label{p3}
Let $(M^{n,r},g)$ be an oriented  (pseudo)-Riemannian manifold, where $g$
is a metric of signature $(r,n-r)$. The metric $g$ induces a
conformal
structure $c:=[g]$ on $M^{n,r}$, which is by definition the  equivalence
class of metrics, which differ from $g$ only by 
multiplication with a positive function in $C^\infty(M)$. 
Such a conformal structure on $M$ is equivalently defined by 
a reduction of the general linear frame bundle $GL(M)$ to a 
principal fibre bundle $G_0(M)$ with structure group
$\mathrm{CO}(r,s)=\RR^+\times\SO(r,s)$. The canonical form with values in
$\RR^{r,s}\cong\frak{g}_-$ reduced to
$G_0(M)$ 
is denoted by $\theta_-$. Moreover, the metric $g$ gives rise to the 
Levi-Civita connection form $\omega_{LC}$ on $G_0(M)$.

The conformal structure $c=[g]$ on $M$ is 
equivalently defined by a 
$P$-reduction $P(M)$ of the second order frame bundle
$GL^{(2)}(M)$, where the structure group $P$ is the parabolic
subgroup of the M{\"o}bius group $\SO(r+1,s+1)$ with
Lie algebra \[\frak{p}:=\frak{g}_0\oplus\frak{g}_1\ .\] The
principal fibre bundle $P(M)$ inherits an invariant 
canonical form $\theta=\theta_-+\theta_0$
from $GL^{(2)}(M)$. Thereby, it is 
\[d\theta_-=[\theta_-,\theta_0]\ ,\]
i.e., the canonical form has no torsion (cf. \cite{Kob72}, \cite{CSS97}). 

Now let $\omega$ be an arbitrary  Cartan connection on $P(M)$ with values 
in
$\frak{so}(r+1,s+1)=\frak{g}_{-1}\oplus\frak{g}_0\oplus\frak{g}_{1}$.
Its curvature $2$-form is defined by 
\[\Omega=d\omega+\frac{1}{2}[\omega,\omega]\]
and the corresponding curvature function with values in
$\frak{g}_-^*\otimes
\frak{g}_-^*\otimes\frak{g}$ is given in a point $u\in
P(M)$ by
\[\nu(u)(x,y):=\Omega(\omega^{-1}_{NC}(x),\omega_{NC}^{-1}(y))(u),\qquad
x,y\in\frak{g}_-\ .\]
Since the map
$ad:\frak{g}_0\rightarrow\frak{gl}(\frak{g}_-)$ is injective,
the $\frak{g}_0$-part $\nu_0$ of the curvature
function can be seen as function on $P(M)$ with values
in $\frak{g}_-^*\otimes
\frak{g}_-^*\otimes\frak{g}_-^*\otimes\frak{g}_+^*$.

It is a well-known fact in conformal geometry
that there exists a unique Cartan connection 
\[\omega_{NC}=\omega_{-1}\oplus\omega_{0}\oplus\omega_{1}\]
on $P(M)$ with the following properties (cf. \cite{Kob72}, \cite{CSS97}):
\begin{enumerate}
\item It is \[\omega_{-1}=\theta_{-1}\qquad\mbox{and}
\qquad\omega_0=\theta_0\ ,
\]
i.e., the torsion of $\omega_{NC}$ vanishes and
\item
\[tr(\nu_0)(x,y):=\sum_{i=1}^{n} \nu_0(e_i,x)(y)(e_i^\flat)=0\ ,\quad
x,y\in \frak{g}_-,\]
i.e., the trace of the $\frak{g}_0$-part of the curvature function is 
trivial.
\end{enumerate}
The so defined Cartan connection $\omega_{NC}$ on the
reduced bundle $P(M)$ is called the canonical 
normal Cartan connection of conformal geometry and is the basic
object for all considerations in this paper.

We want to describe the normal conformal Cartan connection $\omega_{NC}$
in terms of the metric $g$ in
the conformal class $c$. First, we notice that if 
$\pi:P(M)\to G_0(M)$ denotes the natural projection
then $\theta_-$ projects to the canonical form
on $G_0(M)\subset GL^{(1)}(M)$. Furthermore, the $G_0$-equivariant
lifts $\sigma$ of $G_0(M)$ to $P(M)$ correspond bijectively
to the Weyl connections $\omega^{\sigma}$ on $G_0(M)$ by
\[\omega^\sigma=\sigma^*\theta_0.\]
In particular, if $\sigma^{g}$ is the equivariant lift induced by the 
Levi-Civita connection $\omega_{LC}^g$ then 
the $\frak{g}_0$-part of $\omega_{NC}$ is related to $\omega_{LC}^{g}$
by $\omega_{LC}^g=\sigma^*\theta_0$. It remains
to determine the
$\frak{g}_{1}$-part of $\omega_{NC}$ with respect to $g$. This part
must be 
calculated
from the trace-free condition on the curvature function $\nu_0$ and the 
result is 
\[\omega_1=-\Gamma\circ\theta_{-1}\ ,\]
where the function $\Gamma:P(M)\to\frak{g}^*_{-1}\otimes\frak{g}_1$ is 
the pullback of the so-called $\varrho$-tensor on $(M,g)$, which
is given by the expression  
\[K_g=\frac{1}{n-2}\big(\frac{scal_g}{2(n-1)}-Ric_g\big)\ .\]
Thereby,  $Ric_g$ denotes the
Ricci-tenor
and $scal_g$ is the scalar curvature of $g$.
In short, we see that $\omega_{NC}$ is given with respect 
to $g\in c$ by $\theta_-$, $\omega^g_{LC}$ and $K_g$. This
description is invariant in the sense that for all 
$\tilde{g}$ in $c$ the connection $\omega_{NC}$ is determined 
by these data in the same way. This can explain the importance of
the $\varrho$-tensor in conformal geometry. It transforms
naturally in the conformal class.

However, there is a basic construction, which assigns
to every Cartan 
connection on a principal fibre bundle, 
a usual principal fibre bundle connection
through extension. In our case of conformal Cartan geometry,
this can be done as follows. Let 
\[\frak{M}(M)=P(M)\times_{P}\SO(r+1,s+1)\] 
be the
extended bundle with structure group $\SO(r+1,s+1)$. We
call this bundle the M{\"o}bius frame bundle. With respect
to a metric $g$ and the inclusion of $\frak{so}(r,s)$ in
$\frak{so}(r+1,s+1)$ as described above, we can express this bundle
also as
\[\frak{M}(M)=SO(M,g)\times_{\SO(r,s)}\SO(r+1,s+1)\ ,\]
where $SO(M,g)$ is the orthonormal frame bundle to $g$ on $M$.
Then a local frame $s=(s_1,\ldots,s_n)$ on $M$, which is a
local section in $SO(M,g)$,
has a natural extension to a section $s_c=(s_-,s_+,s_1,\ldots,s_n)$
in $\frak{M}(M)$. Thereby, the $(s_{c})_i$'s can be seen as 
tractors (sections) in $T_{\frak{M}}(M)\cong\Lambda^1_{\frak{M}}(M)$ (see 
below). 

The Cartan connection $\omega_{NC}$ can now be extended to
a usual principal connection on $\frak{M}(M)$ by right translation
on the fibres.  
We denote this normal conformal
connection on $\frak{M}(M)$ also by $\omega_{NC}=\omega_{-1}\oplus
\omega_{0}\oplus\omega_1$. 
We have already 
calculated the components of the connection $\omega_{NC}$ with
respect to the metric $g$. Let $s$ be a local frame on $(M,g)$.
Then we have the following expression for the local connection
form on $(M,g)$:
\[ \omega_{NC}\circ ds_c(X)=e_-^\flat\wedge 
\theta_-(X)^\flat+\omega_{LC}\circ
ds_c(X)-e_{+}^\flat\wedge \theta_-(K_g(X))^\flat\ ,\quad X\in TM\ 
,\]
where $(e_-,e_+,e_1,\ldots,e_n)$ is the standard basis in $\RR^{r+1,s+1}$
and $\theta_-$ is evaluated at $s$. The reason for using
the extended approach for $\omega_{NC}$ is because, 
in the following, we would like 
to use
the usual notion of holonomy for a connection on a
principal fibre bundle. We denote the holonomy of $\omega_{NC}$
on $\frak{M}(M)$ over the conformal space
$(M,c)$ with $Hol(\omega_{NC},c)$ (or just $Hol(\omega_{NC})$ if
there is no ambiguity).

Moreover, with the approach of principal connection forms we can introduce 
covariant derivatives to $\omega_{NC}$
on vector bundles with structure group $\SO(r+1,s+1)$ 
associated to $\frak{M}(M)$ in the usual manner.
In particular, $\omega_{NC}$ induces derivatives $\nabla^{NC}$ on the
M{\"o}bius $p$-form bundles (tractor bundles) defined as
\[\Lambda^p_{\frak{M}}(M):=\frak{M}(M)\times_\iota\Lambda^p_{r+1,s+1}\
.\]
With respect to the metric $g$ these bundles split into sums 
of the usual $p$-form bundles on $M$: 
\[\Lambda^{p+1}_{\frak{M}}(M)=\Lambda^p(M)\oplus\Lambda^{p+1}(M)\oplus
\Lambda^{p-1}(M)\oplus\Lambda^p(M)\ .\]
The covariant derivative $\nabla^{NC}$ acts on
sections in these bundles with respect to the above
splitting 
by the matrix expression
\[\nabla^{NC}_X\alpha=
\left(\begin{array}{cccc}\nabla^{LC}_X&-X\ecke&X^\flat\wedge&0\\[2mm]
-K(X)^\flat\wedge&\nabla^{LC}_X&0&X^\flat\wedge\\[2mm]
\ \:\: K(X)\ecke&0&\nabla^{LC}_X&X\ecke\\[2mm]
0&K(X)\ecke&K(X)^\flat\wedge&\nabla^{LC}_X\end{array}\right)
\left(\begin{array}{c}\alpha_-\\[2mm]\alpha_0\\[2mm]\alpha_{\mp}\\[2mm]
\alpha_+
\end{array}\right)\ . \]
Thereby, $\nabla^{LC}$ denotes the Levi-Civita connection. This
expression is calculated straightforwardly  from the local form of
$\omega_{NC}$ and the formulae for the action of $\frak{so}(r+1,s+1)$
on $\Lambda_{r+1,s+1}^{p+1}$. 

%%%%%%%%%%%%%%%%%%%%%%%%%% 
\section{The twistor equations}%%%
%%%%%%%%%%%%%%%%%%%%%%%%%%
\label{p4} 
Let $(M^{n,r},g)$ be an oriented (pseudo)-Riemannian manifold 
and let
$\Lambda^{p+1}_{\frak{M}}(M)$ 
be the associated bundle of $(p+1)$-forms
to
the principal fibre
bundle $\frak{M}(M)$ 
with normal conformal covariant derivative
$\nabla^{NC}$.
We call a section $\alpha\in\Omega_{\frak{M}}^{p+1}(M)$ a 
(normal) twistor
iff
$\nabla^{NC}\alpha=0$. The twistor $\alpha$ corresponds 
via the metric $g$ to
a set
of differential forms on $M^{n,r}$:
\[\alpha\quad\longleftrightarrow\quad 
(\alpha_-,\alpha_0,\alpha_{\mp},\alpha_+)\ ,\]
where $\alpha_-,\alpha_+\in\Omega^p(M)$, $\alpha_0\in\Omega^{p+1}(M)$
and $\alpha_{\mp}\in\Omega^{p-1}(M)$.
The condition $\nabla^{NC}\alpha=0$
is then equivalent to the set of 
conformally covariant equations given by
\begin{eqnarray}
\nabla_X^{LC}\alpha_--X\ecke\alpha_0+X^\flat\wedge\alpha_{\mp}&=&0
\label{g1}\\[2mm]
-K(X)^\flat\wedge\alpha_-+\nabla_X^{LC}\alpha_0
+X^\flat\wedge\alpha_+&=&0\label{g2}\\[2mm]
K(X)\ecke\alpha_-+\nabla_X^{LC}\alpha_{\mp}+X\ecke\alpha_+&=&0
\label{g3}\\[2mm]
K(X)\ecke\alpha_0+K(X)^\flat\wedge\alpha_{\mp}+\nabla_X^{LC}\alpha_+&=&0\
.\label{g4}
\end{eqnarray}          
%For a twistor $2$-form $\alpha$, equation (\ref{g1}) means that the
%$1$-form
%$\alpha_-$ is dual to a conformal vector field. In general, 
%solutions of (\ref{g1}) are known as conformal 
%Killing $p$-forms (cf. \cite{}).
%The equations (\ref{g2}) - (\ref{g4}) impose further 
%conditions on a 
%conformal Killing $p$-form.
%These equations are conformally covariant.
 
We calculate from $\alpha_-$ of a given solution $\alpha$
the remaining differential
forms in order to get equations for $\alpha_-$ only. 
It is 
\[d=\sum_{i=1}^n\varepsilon_i
s_i^\flat\wedge\nabla^{LC}_{s_i}\qquad\mbox{and}\qquad
d^*=-\sum_{i=1}^n\varepsilon_i
s_i\ecke\nabla^{LC}_{s_i}\] 
the exterior differential resp. the codifferential
with respect to a local
orthonormal frame $s$. 
The equations (\ref{g1}) -
(\ref{g3}) imply
for a twistor $\alpha$ of degree $p+1$ that
\begin{eqnarray*}
d\alpha_-&=&(p+1)\alpha_0,\qquad\quad
d^*\alpha_-=(n-p+1)\alpha_{\mp}\\[3mm]
\frac{1}{p+1}d^*d\alpha_-&=&(n-p)\alpha_+-\sum_i^n\varepsilon_i
s_i\ecke(K(s_i)^\flat\wedge\alpha_-)\\[1mm]
\frac{1}{n-p+1}dd^*\alpha_-&=&-p\alpha_+-\sum_i^n\varepsilon_i
s_i^\flat\wedge(K(s_i)\ecke\alpha_-)\ .
 \end{eqnarray*}
For $n\neq 2p$ the sum of the latter two equations results to
\[ 
\alpha_+=\frac{1}{n-2p}\cdot\big(-\frac{scal}{2(n-1)}\alpha_-+
\frac{1}{p+1}d^*d\alpha_-+\frac{1}{n-p+1}dd^*\alpha_-\big), \]
which is 
\[\alpha_+=\frac{1}{n-2p}(\nabla^*\nabla-\frac{scal}{2(n-1)})\alpha_-
\ ,\]
where $\nabla^*\nabla$ denotes the Bochner-Laplacian.
In the middle dimensional case $n=2p$
we have 
\[\alpha_+=\frac{1}{n}\cdot\big[
\frac{1}{p+1}
(d^*d-dd^*)\alpha_-+
\!\sum_{i=1}^n\!\varepsilon_i\cdot\big(s_i\ecke\!(K(s_i)^\flat\wedge
\alpha_-)- 
s_i^\flat\wedge(K(s_i)\ecke\!\alpha_-)\big)\big]\ .\]
We observe that $\alpha_-\equiv 0$ if and only if the twistor
$\alpha$ is trivial.
  
With the so derived expressions for
the components of a twistor $\alpha$ we now formulate 
the
normal twistor equations for a $p$-form
$\alpha_-$
on a (pseudo)-Riemannian manifold $(M^{n,r},g)$. They are
\begin{eqnarray}
0&=&\nabla_X^{LC}\alpha_--\frac{1}{p+1}X\ecke d\alpha_-+ 
\frac{1}{n-p+1}X^\flat\wedge d^*\alpha_{-}
\label{Gtw1}\\[2mm]
0&=&-K(X)^\flat\wedge\alpha_-+\frac{1}{p+1}\nabla_X^{LC}d\alpha_-
+X^\flat\wedge\Box_{p}\alpha_-
\label{Gtw2}\\[2mm]
0&=&K(X)\ecke\alpha_-+\frac{1}{n-p+1}\nabla_X^{LC}d^*\alpha_{-}
+X\ecke\Box_{p}\alpha_-
\label{Gtw3}\\[2mm]
0&=&\frac{1}{p+1}K(X)\ecke d\alpha_-+\frac{1}{n-p+1}
K(X)^\flat\wedge
d^*\alpha_-+\nabla_X^{LC}\Box_{p}\alpha_-\ ,
\label{Gtw4}
\end{eqnarray}  
whereby we set
\[\Box_p:=\frac{1}{n-2p}\cdot
\big(-\frac{scal}{2(n-1)}\mbox{id}+\nabla^*\nabla\big)  
\qquad\mbox{for}\ n\neq 2p
\]
and
\[\Box_{n/2}:=
\frac{1}{n}\cdot\big[
\frac{1}{p+1}
(d^*d-dd^*)+
\sum_{i=1}^n\varepsilon_i\cdot\big(s_i\ecke(K(s_i)^\flat\wedge
\cdot)-
s_i^\flat\wedge(K(s_i)\ecke\cdot)\big)\big]\ .\]
In the following, we say that a $p$-form
$\alpha_-\in\Omega^p(M)$,
which satisfies the (normal) twistor equations (\ref{Gtw1})
- (\ref{Gtw4}), is a normal
conformal Killing $p$-form (or shortly a nc-Killing
$p$-form).
The conformal covariance of the equations implies
that
if $\alpha_-$ is a nc-Killing $p$-form to $g$ on $M$
then the rescaled $p$-form
\[\tilde{\alpha}_-:=e^{-(p+1)\phi}\cdot\alpha_-\]
is 
nc-Killing with respect to the conformally changed metric
$\tilde{g}=e^{-2\phi}\cdot g$. 

However, the
equations (\ref{Gtw1}) - (\ref{Gtw4}) are not only conformally
covariant, but 
a further natural symmetry appears.
Let $*$ denote the Hodge-star 
operator on $\Lambda^*(M)$ defined by
\[\alpha_-\wedge*\alpha_-=g(\alpha_-,\alpha_-)dM\ ,\]
where $dM$ denotes the volume form of $(M^{n,r},g)$. It is
\[**|_{\Lambda^p}=(-1)^{p(n-p)+r}\qquad\mbox{and}\qquad
d^*=(-1)^{n(p-1)+r+1}*d*\ .\]
There is also a `Hodge' operator $*_{\frak{M}}$ on   
$\Lambda_{\frak{M}}^*(M)$
defined
in the same
manner:
\[\alpha\wedge*_{\frak{M}}\alpha=c_{\frak{M}}(\alpha,\alpha)
dM_{\frak{M}}\ ,\]
where $dM_{\frak{M}}:=-e_-^\flat\wedge e_+^\flat\wedge dM$ and
$c_{\frak{M}}$
is the obvious $\SO(r+1,s+1)$-invariant scalar product on
$\Lambda_{\frak{M}}^*(M)$.
The operator $*_{\frak{M}}$ is parallel:
\[\nabla^{NC}*_{\frak{M}}=*_{\frak{M}}\nabla^{NC}\ .\]
Therefore,
if $\alpha$ is a $(p+1)$-twistor then   
$*_{\frak{M}}\alpha$ is a $(n-p+1)$-twistor. The twistor
 $*_{\frak{M}}\alpha$ corresponds to the set
\[  (\ (-1)^p*\alpha_-\ ,\ *\alpha_{\mp}\ ,\
-\!*\alpha_0\ ,\
(-1)^{p+1}*\alpha_+\ )\]
of differential forms. This shows that 
if $\alpha_-$ is a nc-Killing $p$-form
then $*\alpha_-$ is a nc-Killing $(n-p)$-form.
Indeed, 
with
\[
*(X\ecke\beta^p)=(-1)^{p+1}X^\flat\wedge*\beta\qquad 
\mbox{and}\qquad *(X^\flat\wedge\beta^p)=(-1)^pX\ecke
*\beta\ ,\]
and since 
$*\Box_p=-\Box_{n-p}*$ is anti-commuting, 
the normal
twistor equations (\ref{Gtw1}) - (\ref{Gtw4}) are 
seen to be $*$-invariant as well.  

Finally, we remark that for a $1$-form $\alpha_-$,
equation (\ref{Gtw1}) just means that the dual
to $\alpha_-$ is a conformal vector field. In general,
solutions of (\ref{Gtw1}) are known as conformal Killing
$p$-forms (cf. \cite{Kas68}, \cite{Tac69}, \cite{Sem01}). Equation
(\ref{Gtw1}) is Hodge
$*$-invariant itself.
The additional equations (\ref{Gtw2}) -
(\ref{Gtw4}) then impose further conditions on a conformal
Killing $p$-form to be `normal'.

%In particular, this shows that a nc-Killing function
%$\alpha_-$
%without zeros can be rescaled to the constant 
%nc-Killing function
%$\tilde{\alpha}_-=1$
%with respect to $\tilde{g}=\frac{1}{\alpha_-^2}\cdot g$.
%The twistor equations for $\tilde{\alpha}_-$
%immediately show that $\tilde{g}$ is an Einstein metric.

%%%%%%%%%%%%%%%%%%%%%%%%%%%%%%%%%%%%%
\section{Curvature conditions}%%%
%%%%%%%%%%%%%%%%%%%%%%%%%%%%%%%%%%%%%
\label{p5}
We derive here integrability conditions for the existence 
of nc-Killing $p$-forms on a (pseudo)-Riemannian manifold
$(M^{n,r},g)$ in terms of curvature expressions. We denote by 
$s=(s_1,\ldots,s_n)$ a local frame. Let 
\[R(X,Y)Z=\nabla_X\nabla_YZ-\nabla_Y\nabla_XZ-\nabla_{[X,Y]}Z\]
be the Riemannian curvature tensor, where $X,Y,Z\in TM$ are tangent
vectors. By contraction, we obtain
\[ Ric(X)=\sum_{i=1}^n\varepsilon_i\cdot
R(X,s_i)s_i,\qquad scal=tr(Ric)\
,\]
the Ricci tensor and the
scalar curvature. 
The $\varrho$-tensor
is $K=\frac{1}{n-2}\big(\frac{scal}{2(n-1)}-Ric\big)$. The trace-free part
of the
Riemannian curvature tensor is the Weyl tensor $W$, which can be expressed
by 
\[W=R-g\star K\ ,\]
where $\star$ denotes the Kulkarni-Nomizu product.
Moreover, we have the Cotton-York tensor $C$, which is the
anti-symmetrisation of the covariant derivative of the $\varrho$-tensor:
\[C(X,Y):=(\nabla_XK)(Y)-(\nabla_YK)(X)\ .\]
Furthermore, we find the Bach tensor
\[B(X,Y)=\sum_{i=1}^n\varepsilon_i\cdot \nabla_{s_i}C(X,Y,s_i)
-\sum_{i=1}^n\varepsilon_i\cdot W(K(s_i),X,Y,s_i), \]
where $C(X,Y,Z):=C_X(Y,Z)=g(C(Y,Z),X)$.
The Weyl tensor considered as a symmetric map on the space of $2$-forms
is conformally invariant. The Bach tensor is 
symmetric and divergence-free. 
Moreover, we
have the Bianchi identities
\[\begin{array}{l}R(X,Y)Z+ R(Y,Z)X+ R(Z,X)Y=0,\\[3mm]
 \nabla_XR(Y,Z)+ \nabla_YR(Z,X)+ \nabla_ZR(X,Y)=0
\end{array}\]
for all $X,Y,Z\in TM$, which also imply
\[
\sum_{i=1}^n\varepsilon_i\cdot \nabla_{s_i}W(X,Y,Z,s_i)=(3-n)\cdot
C(Z,X,Y)
\quad\ \ \mbox{and}\qquad
\sum_{i=1}^n\varepsilon_i\cdot C(s_i,s_i,X)=0\ .\]
The Kulkarni-Nomizu product $g\star K$ acts on $2$-forms by
\[g\star K(s_i^\flat\wedge s_j^\flat)=s_i^\flat\wedge K(s_j)^\flat-
s_j^\flat\wedge K(s_i)^\flat \ .\]

We calculate now the curvature of the normal conformal
covariant derivative $\nabla^{NC}$ on $\Lambda^{p+1}_{\frak{M}}(M)$. 
For this, let $\alpha=(\alpha_-,\alpha_0,\alpha_{\mp},\alpha_+)$
be a smooth section in $\Lambda^{p+1}_{\frak{M}}(M)$. It is
\begin{eqnarray*}(\nabla^{NC}_X\nabla^{NC}_Y\alpha)_-&=&
\nabla^{LC}_X(\nabla^{LC}_Y\alpha_--Y\ecke\alpha_0+Y^\flat\wedge\alpha_{\mp})
\\[1mm]&&
-X\ecke(-K(Y)^\flat\wedge\alpha_-+\nabla^{LC}_Y\alpha_0+Y^\flat\wedge\alpha_+)
\\[1mm]&&
+X^\flat\wedge(K(Y)\ecke\alpha_-+\nabla_Y\alpha_{\mp}+Y\ecke\alpha_+)\\[3mm]
(\nabla^{NC}_X\nabla^{NC}_Y\alpha)_0&=&
-K(X)^\flat\wedge(\nabla^{LC}_Y\alpha_--Y\ecke\alpha_0+
Y^\flat\wedge\alpha_{\mp})\\[1mm]&&
+\nabla_X^{LC}(-K(Y)^\flat\wedge\alpha_-+
\nabla^{LC}_Y\alpha_0+Y^\flat\wedge\alpha_+)\\[1mm]&&
+X^\flat\wedge(K(Y)\ecke\alpha_0+K(Y)^\flat\wedge\alpha_{\mp}+
\nabla_Y^{LC}\alpha_+)\\[3mm]
(\nabla^{NC}_X\nabla^{NC}_Y\alpha)_{\mp}&=&
K(X)\ecke(\nabla^{LC}_Y\alpha_--Y\ecke\alpha_0+
Y^\flat\wedge\alpha_{\mp})\\[1mm]&&
+\nabla_X^{LC}(K(Y)\ecke\alpha_-+\nabla_Y\alpha_{\mp}+Y\ecke\alpha_+)
\\[1mm]&&
+X\ecke(K(Y)\ecke\alpha_0+K(Y)^\flat\wedge\alpha_{\mp}+
\nabla_Y^{LC}\alpha_+)\\[3mm]
(\nabla^{NC}_X\nabla^{NC}_Y\alpha)_+&=&
K(X)\ecke(-K(Y)^\flat\wedge\alpha_-+
\nabla^{LC}_Y\alpha_0+Y^\flat\wedge\alpha_+)\\[1mm]&&
+K(X)^\flat\wedge(K(Y)\ecke\alpha_-+\nabla_Y\alpha_{\mp}+Y\ecke\alpha_+)
\\[1mm]&&
+\nabla^{LC}_X(K(Y)\ecke\alpha_0+K(Y)^\flat\wedge\alpha_{\mp}+
\nabla_Y^{LC}\alpha_+)\end{eqnarray*}
and we obtain 
\begin{eqnarray*}
(R^{\nabla}(X,Y)\circ\alpha)_-&=&R^{LC}(X,Y)\circ\alpha_-
\\[1mm]&&
+(X\ecke(K(Y)^\flat\wedge\alpha_-)-
Y\ecke(K(X)^\flat\wedge\alpha_-))\\[1mm]&&
+(X^\flat\wedge(K(Y)\ecke\alpha_-)
-Y^\flat\wedge(K(X)\ecke\alpha_-))\\[1mm]
&=& W(X,Y)\circ\alpha_-\\[3mm]
(R^{\nabla}(X,Y)\circ\alpha)_0&=&
W(X,Y)\circ\alpha_0-C(X,Y)^\flat\wedge\alpha_-\\[3mm]
(R^{\nabla}(X,Y)\circ\alpha)_{\mp}&=&
W(X,Y)\circ\alpha_{\mp}+C(X,Y)\ecke\alpha_-\\[3mm]
(R^{\nabla}(X,Y)\circ\alpha)_+&=&
W(X,Y)\circ\alpha_++C(X,Y)\ecke\alpha_0+C(X,Y)^\flat\wedge\alpha_{\mp}\ ,
\end{eqnarray*}
i.e., the curvature takes the matrix form
\[R^\nabla=
\left(\begin{array}{cccc}W&0&0&0\\[2mm]
-C(X,Y)^\flat\wedge&W&0&0\\[2mm]
\ \:\: C(X,Y)\ecke&0&W&0\\[2mm]
0&C(X,Y)\ecke&C(X,Y)^\flat\wedge&W\end{array}\right)
\ . \]
As integrability condition for the existence of a twistor $\alpha$
we obtain
\begin{eqnarray*}
W(X,Y)\circ\alpha_-&=&0\\[3mm]
W(X,Y)\circ\alpha_0&=&\ \:C(X,Y)^\flat\wedge\alpha_-\\[3mm]
W(X,Y)\circ\alpha_{\mp}&=&-C(X,Y)\ecke\alpha_-\\[3mm]
W(X,Y)\circ\alpha_+&=&-C(X,Y)\ecke\alpha_0-C(X,Y)^\flat\wedge\alpha_{\mp}\
.
\end{eqnarray*}
By taking the divergence on both sides 
of these equations we get
\begin{eqnarray*}
(n-4)\cdot C_T\circ\alpha_-&=&0\\[3mm]
(n-4)\cdot C_T\circ\alpha_0&=&-B(T)^\flat\wedge\alpha_-\\[3mm]
(n-4)\cdot C_T\circ\alpha_{\mp}&=&\ \: B(T)\ecke\alpha_-\\[3mm]
(n-4)\cdot C_T\circ\alpha_+&=&\
\: B(T)\ecke\alpha_0+B(T)^\flat\wedge\alpha_{\mp}\
.
\end{eqnarray*}
The curvature conditions for 
the nc-Killing $p$-form $\alpha_-$ take the form:
\begin{eqnarray}
\!\!\!\!\!\!W(X,Y)\circ\alpha_-&\!\!\!=\!\!&\!0 \label{G17}\\[3mm]
\!\!\!\!\!\!W(X,Y)\circ d\alpha_-&\!\!\!=\!\!&\!\ \
(p+1)\cdot\label{G18}
C(X,Y)^\flat\wedge\alpha_-\\[3mm]   
\!\!\!\!\!\!W(X,Y)\circ d^*\alpha_-&\!\!\!=\!\!&\!-(n-p+1)C(X,Y)
\ecke\alpha_-\label{G19}\\[3mm]
\!\!\!\!\!\!W(X,Y)\circ\Box_p\alpha_-&\!\!\!=\!\!&\!
-\frac{1}{p+1}
C(X,Y)\ecke d\alpha_--\frac{1}{n-p+1}C(X,Y)^\flat\wedge
d^*\alpha_{-}. \label{G20}
\end{eqnarray}
Of course, the sets of integrability 
conditions are conformally covariant and
invariant under the Hodge $*$-operator.\\

%\[ C_{\tilde{g}}=C_g+(n-3)grad_g\phi\ecke W\ .\]

%%%%%%%%%%%%%%%%%%%%%%%%%%%%%%%%%%%%%%%%%%%%%%%%%%%%%%%%%%%%%
\section{Normal conformal Killing $p$-forms on Einstein manifolds}%%%
%%%%%%%%%%%%%%%%%%%%%%%%%%%%%%%%%%%%%%%%%%%%%%%%%%%%%%%%%%%%%
\label{p6}
We consider in this paragraph solutions of the 
normal twistor equations (\ref{Gtw1}) - (\ref{Gtw4})
on Einstein manifolds.
Before we start with this, we want to state a criterion when
a space $(M^{n,r},g)$ is conformally equivalent to
an Einstein space, i.e.,  
there is a
metric $\tilde{g}$ in the conformal class $c=[g]$,
which satisfies 
\[Ric_{\tilde{g}}=\frac{scal_{\tilde{g}}}{n}\cdot
\tilde{g}\ .\]

For this, let us
assume that $f_-=\alpha_-$ is a nc-Killing function (= $0$-form)
without
zeros. We have mentioned before that the rescaled function
$\tilde{\alpha}_-=\frac{1}{f_-}\alpha_-=1$ is nc-Killing
with respect to the metric 
$\tilde{g}=\frac{1}{f^2_-}\cdot g$.
From the twistor equations (\ref{Gtw1}) - (\ref{Gtw4}), 
it follows immediately 
\[K_{\tilde{g}}=-\frac{scal_{\tilde{g}}}{2n(n-1)}
\cdot\tilde{g}\ ,\] 
which means that $\tilde{g}$ is
Einstein. On the other hand, every constant
function on an Einstein space is nc-Killing. The criterion then says that 
a metric is conformally Einstein (i.e., there is 
an Einstein metric in the conformal class $[g]$) if and only if there
exists at least one nc-Killing function without zeros. 

Obviously, in case that $f_-$ has a zero the rescaling
in the way as above is not possible.
Indeed, examples of nc-Killing functions 
on non-conformally Einstein spaces 
%(that means there is no Einstein metric
%in the conformal class) 
are well known (cf. 
\cite{KR96}). However, since in general the set of zeros of nc-Killing
forms
is singular on the underlying conformal space, we can at least say
that the existence of a nc-Killing function $f_-$
implies that, up to singularities, 
an Einstein metric exists in the conformal class $[g]$.
This is exactly the case when the holonomy group of the conformal
connection $\omega_{NC}$
fixes at least one vector in $\RR^{r+1,s+1}$.
In case that this vector is lightlike the `Einstein scaling' is 
Ricci-flat.
The timelike case is for $scal>0$, the spacelike when $scal<0$.
By the way, 
the  normal twistor equations 
for a function $f_-$ are in general equivalent to
\[\begin{array}{l}Hess(f_-)=f_-\cdot K_o\qquad\mbox{and}
\qquad K(X)(f_-)=\frac{1}{n}X(\frac{scal}{2(n-1)}f_-)
\end{array}\] 
for all $X\in TM$, 
where $K_o$ denotes the trace-free part
of the $\varrho$-tensor.

% and $\Delta=d^*d$
%is the usual Laplacian on functions. 
%If a solution
%$f_-$ of these equations has a zero on an Einstein
%space then this space has constant sectional 
%curvature (cf. \cite{}).
%Next we want to mention a further case when
%the existence of a nc-Killing form implies
%the Einstein condition (more exactly here the Ricci-flatness).
%For this, let $M^{n,r}$ be a space of even dimension
%($n=2m$) and signature. In this case the bundle 
%$\Lambda^m(M)$ of forms of degree
%$m$ in the middle dimension splits into a selfdual and
%an anti-selfdual part: $\Lambda^m_+(M)\oplus
%\Lambda^m_-(M)$. It is known
%that an (anti)-selfdual conformal Killing $m$-form 
%$\alpha_-$, which has no zeros and
%whose length $\phi=\sqrt{\pm g(\alpha_-,\alpha_-)}$ does 
%not vanish, is conformally equivalent to the parallel form
%$\tilde{\alpha}_-=e^{-(m+1)\phi}\alpha_-$  
%with respect to the metric $\tilde{g}=e^{-2\phi}g$ 
%(cf. \cite{}).
%In general, as one can 
%immediately see from the twistor equations (\ref{Gtw2}) -
%(\ref{Gtw3}),
%the existence of a 
%parallel nc-Killing form of degree $p\neq 0,n$ 
%with non-vanishing length
%implies
%that the $\varrho$-tensor $K$ vanishes,   
%i.e., the underlying metric  is Ricci-flat. In 
%case of the existence of 
%an (anti)-selfdual nc-Killing form $\alpha_-$
%with non-vanishing
%length and without zeros this fact leads to the conclusion
%that the underlying metric $g$
%is conformally equivalent to a Ricci-flat metric.    

Now we assume that 
$(M^{n,r},g)$ is a (pseudo)-Riemannian Einstein manifold. 
The constant functions
are nc-Killing on Einstein manifolds.
In particular, the 1-form
$o:=s_-^\flat-\frac{scal}{2(n-1)n}s_+^\flat$, which 
obviously satisfies the equations (\ref{g1}) - (\ref{g4}), 
is the normal twistor in $\Omega^1_{\frak{M}}(M)$,
which corresponds to the set 
of forms $(1,0,0,-\frac{scal}{2(n-1)n})$, i.e., $o_-=1$
is the constant unit function. 
Furthermore, let
\[\alpha=s_-^\flat\wedge\alpha_-+\alpha_0+s_-^\flat\wedge
s_+^\flat\wedge\alpha_{\mp}+s_+^\flat\wedge\alpha_+\]
be an arbitrary normal twistor of degree $(p+1)$ 
with $d\alpha_-\not\equiv 0$ on $M$.
It follows 
immediately that $o\wedge\alpha$ is a $(p+2)$-twistor. 
This twistor corresponds to the set
\[(\ \alpha_0\ ,\ 0\ ,\ \alpha_++\frac{scal}{2(n-1)n}\alpha_-\
,\
\frac{-scal}{2(n-1)n}\alpha_0\ ) \]
of differential forms,
which shows that $d\alpha_-$ is a (closed)
nc-Killing $(p+1)$-form.  
The $(n-p-1)$-form $*d\alpha_-$
is then nc-Killing 
and coclosed. 

In general, the set of 
twistor equations (\ref{Gtw1}) - (\ref{Gtw4})
reduces for a coclosed $p$-form
$\beta_-$ on
an Einstein space to 
\begin{eqnarray*}\nabla^{LC}_X\beta_-&=&\frac{1}{p+1}\cdot X\ecke
d\beta_-\ ,\\[2mm]
\nabla^{LC}_Xd\beta_-&=&-\frac{(p+1)\cdot 
scal}{n\cdot (n-1)}\cdot X^\flat\wedge\beta_-\ ,
\end{eqnarray*}
which implies $\Delta_p\beta_-=
\frac{(p+1)(n-p)\cdot 
scal}{n\cdot(n-1)}\beta_-$ for the Laplacian 
$\Delta_p=dd^*+d^*d$.
A differential form that satisfies these two equations above
(not only in the Einstein case) is
called
a special Killing $p$-form 
to the Killing constant
$-\frac{(p+1)\cdot scal}{n\cdot (n-1)}$
(cf. \cite{Sem01}). 
There is a nice way to describe the geometry of spaces with
special Killing forms for a non-zero Killing constant using
the cone construction. We explain this approach briefly as next.
In a further step we show that we can extend this approach
to describe all nc-Killing 
forms on Einstein spaces with non-zero scalar curvature 
as parallel forms on an `ambient' metric.  

The cone metric with scaling $b\neq 0$ is defined on the space  
$\RR_+\times M$
as \[\hat{g}_b:=bdt^2+t^2g\ .\] This metric
has either signature $(r,s+1)$ or $(r+1,s)$. 
We have the following result for special Killing forms.
\begin{PR} \label{PR49}(cf. \cite{Sem01}) 
Let  
$(M^{n,r},g)$ be a 
(pseudo)-Riemannian manifold and $\hat{M}_b$ its cone
with metric $\hat{g}_b$ to the constant
$b=\frac{(n-1)n}{scal_g}\neq
0$.  Then 
the special Killing $p$-forms on $M^{n,r}$ with
Killing constant $-\frac{(p+1)\cdot scal}{n\cdot(n-1)}$
correspond bijectively
to the parallel
$(p+1)$-forms on the cone $\hat{M}_b$. 
The correspondence is given by
\[\beta_-\in\Omega^p(M)
\qquad\mapsto\qquad t^pdt\wedge\beta_-
+\frac{sign(b)\cdot t^{p+1}}{p+1}d\beta_-\in\Omega^{p+1}(\hat{M}_b)\ .\] 
\end{PR}
However, the metric $\bar{g}_b$  on $\RR_+^2\times M$ 
of signature 
$(r+1,s+1)$ 
in dimension $n+2$ defined
by
\[\bar{g}_b:=b(dt^2-ds^2)+t^2\cdot g\]
is appropriate to describe all normal twistors 
on Einstein spaces ($Ric\neq 0$) as parallel forms.
\begin{PR} 
\label{PR50} Let $(M^{n,r},g)$ be an Einstein space and $(\bar{M}_b,
\bar{g}_b)$ its ambient metric with $b=\frac{(n-1)n}{scal}\neq 0$.
There is a natural and bijective correspondence between normal
twistors $\alpha\in\Omega_{\frak{M}}^{p+1}(M)$ and parallel forms
$\bar{\alpha}\in\Omega^{p+1}(\bar{M}_b)$. Moreover, the holonomy
groups of the normal conformal connection and the Levi-Civita connection
on $\bar{M}_b$ coincide, i.e.,\[
Hol(\omega_{NC},c)\cong Hol(\omega_{LC}^{\bar{g}_b})\ .\]
\end{PR}

{\bf Proof.} We prove the statement that the 
holonomies coincide.
For this, we embed $M$ in $\bar{M}_b$ by $i:M\to M\times\{(1,1)\}$. 
The $1$-twistor on $M$ which defines the Einstein structure is
$s^\flat_--\frac{scal}{2(n-1)n}s_+^\flat$.
Next we define an isometric map between the tangent tractors on $M$ and 
the tangent vectors
at $i(M)\subset \bar{M}_b$ by assigning
\[\begin{array}{ccr}
s_i\qquad&\to&\qquad s_i\quad \in T\bar{M}_b|_{M\times\{(1,1)\}}\\[3mm]
\sqrt{|b|}\cdot(
s_-+\frac{scal}{2(n-1)n}s_+^\flat)\qquad&\to&\qquad 
\sqrt{\frac{1}{|b|}}\cdot\partial_t\quad \in
T\bar{M}_b|_{M\times\{(1,1)\}}\\[3mm]
\sqrt{|b|}\cdot(
s_--\frac{scal}{2(n-1)n}s_+^\flat)\qquad&\to&\qquad
\sqrt{\frac{1}{|b|}}\cdot\partial_s\quad \in
T\bar{M}_b|_{M\times\{(1,1)\}}\end{array}\]
It can be easily calculated that the resulting bundle isomorphism
between $\frak{M}(M)$ and $SO(\bar{M}_b)|_{M\times\{(1,1)\}}$
has the property $i^*\omega_{LC}^{\bar{g}}=\omega_{NC}$, i.e., the
connection forms are identified in this way.
Moreover, all elements in the holonomy of the metric 
$\bar{g}$
are generated by the horizontal lifts of paths which move solely
on the level set $M\times\{(1,1)\}$. This shows that 
$Hol(\omega_{NC},c)\cong Hol(\omega_{LC}^{\bar{g}_b})$. 

The above map between the tangent tractors on $M$ and the vectors on 
$\bar{M}_b$ also provides the correspondence between the normal
twistors on $M$
and the parallel forms on the ambient space $\bar{M}_b$.\hfill$\Box$\\

In particular, one can see from Proposition \ref{PR50} that every nc-Killing
form on an Einstein space with non-zero scalar curvature is the sum
of a closed and a coclosed nc-Killing form. The holonomy
groups of the Levi-Civita connections on the cone $\hat{M}_b$ and the 
ambient $\bar{M}_b$ are obviously identical, which also means
$Hol(\omega_{NC})=Hol(\omega_{LC}^{\hat{g}_b})$.

In the Riemannian case 
a geometric characterisation
of complete spaces $(M^n,g)$ with positive scalar curvature
admitting special Killing forms was established
by using the above correspondence with the cone and
the holonomy classification  
for simply connected, irreducible and  non-locally symmetric
spaces  
(cf. \cite{Bar93}, \cite{Sem01}).
Thereby, we remember to the fact that if the holonomy of
a Riemannian cone $\hat{M}_b$, $b>0$, over a 
complete Riemannian space $M$ 
is reducible then the cone is automatically flat.
However, it is not difficult to extend the geometric
characterisation to non-complete spaces with special
Killing forms
when the cone is reducible, but not flat.
In this case the metric $g$ on $M$ turns out to be  (locally) 
a certain warped-product. Combining Proposition \ref{PR49}
and \ref{PR50} results to the following.

\begin{THEO}\label{TH1} a)
Let $(M^n,g)$ be a simply connected and complete
Riemannian Einstein space of positive scalar curvature
admitting
a nc-Killing $p$-form. Then $M^n$ is either
\begin{enumerate} 
\item the round (conformally flat) sphere $S^n$,
\item an Einstein-Sasaki manifold of 
odd dimension $n\geq 5$ with 
a special Killing $1$-form $\alpha_-$,
\item an Einstein-3-Sasaki space of 
dimension $n=4m+3\geq 7$ with three
independent special Killing $1$-forms $\alpha_-^1$, $\alpha_-^2$ and
$\alpha_-^3$,
\item a nearly K{\"a}hler manifold of
dimension $6$, where the K{\"a}hler
form $\omega_-$ is a special
Killing $2$-form or
\item a nearly parallel $G_2$-manifold in
dimension $7$ with its fundamental form $\gamma_-$
as special 
Killing $3$-form.
\end{enumerate}
b) If the space $M^n$ is not complete
and the cone reducible then the metric $g$ has up to a constant scaling
factor (locally) the form
\[dt^2+\sin^2(t)\cdot k+\cos^2(t)\cdot h\ ,
 \]
where $k^p$ and $h^q$ are arbitrary Riemannian Einstein
metrics of positive scalar curvature on spaces with dimension
$p$ resp. $q$ ($0\leq q\leq n-1$). 
The scaled volume forms 
\[\sin^{-p}\cdot\, dvol_k\qquad\mbox{and}\qquad
\cos^{-q}\cdot\, dvol_h\]
to $k$ and $h$ are special Killing of degree $p$ resp. $q$. 
\end{THEO}

Similarly, for nc-Killing forms on Riemannian Einstein spaces 
$(M^n,g)$ 
of negative scalar curvature one has to consider
the cone with Lorentzian metric (indefinite of signature 
$(1,n)$). In this case the Lorentzian cone either has
weakly irreducible or decomposable holonomy. 
In both cases one can show that $(M^n,g)$ admits 
(locally) certain warped-product structures (cf. \cite{Bau89}).
In general, there is no
classification of possible holonomy groups 
of the Levi-Civita connection for 
pseudo-Riemannian spaces.
This implicates the lack
of a further geometric characterisation
of pseudo-Riemannian Einstein spaces with nc-Killing forms.

A parallel form on a Ricci-flat metric $g$ is 
nc-Killing. The lightlike $1$-twistor $o=s_-^\flat$ gives 
rise to the constant nc-Killing function $f_-=1$. The holonomy
of the normal conformal connection to $[g]$ will be weakly irreducible,
in 
general.
A reproduction of the normal conformal connection
with its holonomy by some ambient metric with its Levi-Civita connection
and the corresponding holonomy is 
not done here in the general situation.

%In general, a conformal Killing $p$-form $\alpha_-$
%is conformally
%equivalent to a parallel form for some metric 
%$\tilde{g}=e^{-2\phi}\cdot g$ if and only if
%\[ d\alpha_-=(p+1)d\phi\wedge\alpha_-\qquad
%\mbox{and}\qquad d^*=-(n-p+1)\cdot grad(\phi)\wedge\alpha_-
%\ . \]
%For (anti)-selfdual forms these two equations 
%are equivalent. This was used in \cite{} to show
%that any (anti)-selfdual conformal Killing $n/2$-form 
%$\alpha_-$
%with length function $e^{\phi}$ (without zeros)
%on a Riemannian space $M^n$
%is parallel with respect to $\tilde{g}=e^{-2\phi}g$.
%If $\tilde{g}$ is Einstein then
%\[\Box_{n/2}=\frac{1}{n\cdot(n/2+1)}
%(d^*d-dd^*)\ ,\]
%which implies here that 
%$\tilde{\alpha}_+$
%vansihes for (anti)-selfdual $\tilde{\alpha}_-=
%e^{-(n/2+1)\phi}\alpha_-$
%and from the twistor equations (\ref{Gtw2}) and (\ref{Gtw3})
%it follows that $\tilde{g}$ is Ricci-flat.

Finally, we want to take a closer look on the warped-product
structure 
in Theorem \ref{TH1} for the case when the 
Riemannian cone is reducible and has
one parallel vector $\hat{P}$.
Through the correspondence with the cone, 
the vector $\hat{P}$ can be seen to give 
rise to a function $f_{\hat{P}}$ on $(M^n,g)$,
which satisfies the second order differential equation 
\[\nabla df_{\hat{P}}=-\frac{f_{\hat{P}}}{c}\cdot g\ ,\]
i.e., the vector $grad(f_{\hat{P}})$ is a 
conformal gradient field.
It is well-known that the existence
of a conformal gradient field 
gives rise to a (local) 
warped-product structure on 
the Riemannian space $M^n$ (cf. \cite{KR97}).
In Theorem \ref{TH1}, this is the special 
case when the metric
$h$ or $k$ vanishes ($p$ or $q=0$).
Moreover, in case that $f_{\hat{P}}$ or $grad(f_{\hat{P}})$ has a zero on
an
Einstein
space its sectional curvature is constant.
%As we mentioned already zeros of $f_{\hat{P}}$
%may appear at points, where the space is not conformally
%Einstein.

Furthermore, conformal gradient fields are known
to generate conformal transformations between Einstein
spaces (cf. \cite{Bri25}, \cite{Kuh88}). This can be understood with our
approach in
the following way. The function $f_{\hat{P}}$ on an
Einstein space $(M^n,g)$ 
is a nc-Killing function (or special Killing function)
and we mentioned already 
above that the scaling of the metric $g$ by
a nc-Killing function without zeros gives rise to a conformal
transformation to another Einstein metric. 

To summarise,
an Einstein space $(M^n,g)$ has always nc-Killing functions 
(the constant functions). In case that there is in
addition a non-constant nc-Killing function $f_-$ without zeros,
whose corresponding twistor is 
$\alpha\in\Omega^1_{\frak{M}}(M)$,
we find a conformal transformation to a further Einstein
metric $\tilde{g}\in [g]$. 
In particular, the normal $2$-twistor $o\wedge\alpha$
corresponds to the conformal gradient field
$grad(f_-)$, which gives rise to the warped product structure.
More generally, a set of $j$ `independent' nc-Killing functions
on an Einstein space $M^n$ 
induces $j-1$ different conformal transformations to 
further Einstein
metrics in $[g]$ and $j-1$ different ways of
expressing warped-product structures.
In case that $M^n$ is the $n$-sphere $S^n$, there is 
the constant nc-Killing function $o_-$ on $S^n$ and there are 
$n+1$ further `independent' nc-Killing functions,
each of them  with an isolated
zero on $S^n$, which give rise to $n+1$ 
conformal transformations
to Einstein metrics with constant sectional curvature
up to a singularity
(stereographic projections).
This is the conformally flat case, where the number
of `independent' nc-Killing functions 
on a space is the maximal one (i.e., $n+2$).
  
%We have shown in this paragraph that the geometric
%structure of Einstein spaces with nc-Killing forms
%is well known.
%In particular, it turns out that the holonomy
%group of the normal conformal connection fixes
%at least one vector, i.e., the holonomy
%representation is reducible (we will also say decomposable
%in abuse of the common usage).
%By taking a closer look one can see even more. The
%holonomy of $\omega_{NC}$ is related (even determined)
%to the holonomy of the Levi-Civita connection of the 
%Einstein space ($scal=0$) or its cone ($scal\neq 0$).
%In the next paragraph we extend this philosophy
%to the case when there exist decomposable twistors
%of arbitrary degree (not only degree $1$). 

%We want to apply our results in the following to the study
%of nc-Killing forms on (pseudo)-Riemannian spaces
%in dimension $4$. In particular, we attempt to 
%find a holonomy classification for the normal
%conformal connection $\omega_{NC}$. 
%Unsurprisingly, in the Einstein case 
%the resulting holonomy list 
%will be related to the holonomy group of the Levi-Civita
%connection on the cone. However, 
%there are solutions of the normal
%twistor equations for $p$-forms on non-Einstein spaces.
%We will come across such interesting examples 
%when we investigate
%the case of space-times
%in paragraph \ref{}. Fortunately, it seems that
%this effect makes the 
%normal twistor equation and its solutions 
%from the viewpoint of conformal
%geometry more interesting than one could suspect
%from the results of this paragraph.

%%%%%%%%%%%%%%%%%%%%%%%%%%%%%%%%%%%%%%%%%%%%%%%%%%%%%%%%%%%%%%
\section{Solutions with decomposable twistors}%%%
%%%%%%%%%%%%%%%%%%%%%%%%%%%%%%%%%%%%%%%%%%%%%%%%%%%%%%%%%%%%%%
\label{p7}
In this paragraph we want to investigate conformal
spaces $(M^n,[g])$, which admit decomposable normal twistors
\[\alpha=\alpha_1\wedge\cdots\wedge\alpha_{p+1}
\in\Omega^{p+1}_{\frak{M}}(M)\ ,\]
i.e., the $\alpha_i$'s are $1$-forms in $\Omega^{1}_{\frak{M}}(M)$. 
The existence of such a twistor implies that the 
holonomy representation of the normal conformal
connection $\omega_{NC}$ has an invariant 
(non-trivial) subspace
in $\RR^{r+1,s+1}$. That means
the representation is not irreducible. 
We remember that we studied in the previous
paragraph, which was about Einstein metrics, the case of a $1$-twistor.
Of course, a $1$-twistor is always decomposable.

First we want to observe here an easy generation principle for 
coclosed nc-Killing forms from a given one. 
So let 
$\alpha_-$
be such a coclosed nc-Killing form 
on a space $H$ of dimension $p$ with Einstein metric $h$  of scalar 
curvature
$scal_h$. Now we consider the product metric
\[g=h\times l\ ,\]
where $l$ is a metric on a space $L$ of dimension $q$,
and we produce the pullback of $\alpha_-$ to $M=H\times L$.
Obviously, the first of the normal twistor equations (\ref{Gtw1})
for the pullback $\alpha_-$ on $M$ is still satisfied, since 
for every $Y\in TL$ 
\[\nabla_Y\alpha_-=Y\ecke\alpha_0=Y^\flat\wedge\alpha_\mp=0\ .\]
Indeed, it is 
straightforward to show that if we choose $l$ to be Einstein 
with scalar curvature
\[scal_l=-\frac{q(q-1)}{p(p-1)}\cdot scal_h\ ,\]
then the three remaining normal twistor equations (\ref{Gtw2}) -
(\ref{Gtw4}) 
are also satisfied for the pullback and we can conclude that we have 
produced
a metric $g$ on a space of dimension $n=p+q$ admitting
a coclosed nc-Killing form. (If we choose for $l$ a different scalar 
curvature then
$\alpha_-$ is only a conformal Killing form on $h\times l$.
This shows that, of course, not every conformal Killing form 
is normal.) 

Next we formulate and prove a lemma which on the other side shows that
certain conformal Killing forms give rise to a product metric in the 
conformal class of a space.
\begin{LM} \label{Lem11} Let $\alpha_-$ be a conformal Killing $p$-form with 
$\|\alpha_-\|^2\neq 0$
on a space $(M^n,g)$ satisfying the following three properties:
\begin{enumerate}
\item
$\alpha_-$ is decomposable, i.e., $\alpha_-=\alpha_1^1\wedge\ldots\wedge
\alpha_p^1$ is a $\wedge$-product of $p$ $1$-forms,
\item there is  $A\in \Gamma(TM)$ such that 
$d\alpha_-=A^\flat\wedge\alpha_-$
and
\item
there is $B\in \Gamma(TM)$ such that $d^*\alpha_-=B\ecke\alpha_-$.
\end{enumerate} 
Then it exists a rescaled metric $\tilde{g}$ in the 
conformal class
$[g]$ such that the rescaled conformal Killing form $\tilde{\alpha}_-$ is 
parallel. In particular, if $0<p<n$ then $\tilde{g}$ is  
(locally) a product metric $h\times l$.  
\end{LM}
{\bf Proof.} First we observe that the three 
assumptions are invariant under conformal rescaling, e.g. it is
$d(e^\phi\alpha_-)=\tilde{A}^\flat\wedge e^\phi\alpha_-$ with
$\tilde{A}=d\phi+A$.
Moreover, since $\|\alpha_-\|^2\neq 0$, 
we can scale the metric $g$ such that
$\alpha_-$ has constant non-zero length. 
For simplicity, we assume that $g$ is already in this scaling. 
Then it is
\begin{eqnarray*}0&=&X(g(\alpha_-,\alpha_-))=2g(\nabla_X\alpha_-,\alpha_-)\\
&=&2g(\frac{1}{p+1}X\ecke d\alpha_-
-\frac{1}{n-p+1}X^\flat\wedge d^*\alpha_-,\alpha_-)\ .\end{eqnarray*}
But from the assumptions, we see that this is only possible if $A,B=0$, 
i.e., $\alpha_-$ is
closed and coclosed which means that it is parallel,
since it is a conformal Killing form. Moreover, 
$\alpha_-$ is decomposable and this shows that $g$ is (locally)
a product $h\times l$ for the case when $\deg(\alpha_-)\neq 
0,n$.\hfill$\Box$\\

This lemma generalises the well-known fact that a conformal
vector field, which is hypersurface orthogonal, is parallel 
with respect to some metric in the conformal class.
We also remark at this point 
that in general a conformal Killing $p$-form $\alpha_-$
is conformally
equivalent to a parallel form for some metric
$\tilde{g}=e^{-2\phi}\cdot g$ in the conformal class if and only if
\[ d\alpha_-=(p+1)\cdot d\phi\wedge\alpha_-\qquad
\mbox{and}\qquad d^*\alpha_-=-(n-p+1)\cdot grad(\phi)\ecke\alpha_-
\ . \]
This shows that $A,B\neq 0$ in the above lemma are actually 
parallel vectors (if they exist). 
Moreover, for an (anti)-selfdual form the latter two equations 
are equivalent and this can be used to show
that any (anti)-selfdual conformal Killing $(n/2)$-form
$\alpha_-$
with length function $e^{\phi}$ 
on a Riemannian space $M^n$
is parallel with respect to $\tilde{g}=e^{-2\phi}g$ (cf. \cite{Sem01}).

Now we are prepared to consider nc-Killing forms $\alpha_-$ with 
corresponding decomposable twistor. Let $\alpha$ be such a 
twistor of degree $p+1$. 
The first statement that we can make says that 
the four corresponding differential forms  
$\alpha_-$, $\alpha_0$, $\alpha_\mp$ and $\alpha_+$ all are 
decomposable as well. For example, it is
\[\alpha_-=s_+\ecke(s_-\ecke(s^\flat_+\wedge\alpha))\]
and henceforth, $\alpha_-$ is obviously decomposable.
But we can say even more. The twistor $\alpha$ has two 
different normal
forms with respect to a fixed frame 
$(s_-,s_+,s_1,\ldots,s_n)$. In the first case,
it is $s_-\ecke(s_+\ecke\alpha)=0$ and the corresponding normal form is
given by
\[\alpha=(a\cdot s_-^\flat+b\cdot s_+^\flat+c\cdot t^\flat_{p+1})
\wedge t_1^\flat\wedge\ldots \wedge t_p^\flat\ ,\]
where the $t_i$'s are orthogonal to each other 
and are contained in the span of the $s_i$'s, 
$i=1,\ldots, n$, and $a,b$ and 
$c$ are
some constants.  
In the second case, it is $s_-\ecke(s_+\ecke\alpha)\neq 0$
and we have the normal form
\[\alpha=(a\cdot s_-^\flat+b\cdot t_p^\flat+c\cdot t^\flat_{p+1})
\wedge (d\cdot s_+^\flat + t_p^\flat)\wedge t_1^\flat\wedge\ldots 
\wedge t_{p-1}^\flat\ .\] 
From these two normal forms we can see that for a twistor $\alpha$
there are 
vectors $A,B$ such that
\[d\alpha_-=A^\flat\wedge\alpha_-\qquad\mbox{and}\qquad
d^*\alpha_-=B\ecke
\alpha_-\ . \]
Indeed, we can apply now Lemma \ref{Lem11} and obtain the following 
result.
\begin{LM} \label{PR11} Let $\alpha_-$ be a nc-Killing $p$-form 
on $(M,g)$ 
with $\|\alpha_-\|^2\neq 0$ such that
the corresponding normal twistor $\alpha$ is decomposable. Then there
exists $\tilde{g}$ in the conformal class $[g]$ such that 
the rescaled form $\tilde{\alpha}_-$ is parallel. 
\end{LM}
The parallel nc-Killing $p$-form that is guaranteed by Lemma
\ref{PR11}
is decomposable. If its degree is different from $0$ and $n$, it 
gives rise (locally) to a product metric
$h\times l$ in the conformal class $[g]$. Moreover, the normal twistor
equations (\ref{g2}) and (\ref{g3}) show that the factors $h,l$ are
Einstein and we can conclude that $\tilde{g}$ is a product of 
Einstein metrics with 
\[ scal_l=-\frac{(n-p)\cdot(n-p-1)}{p\cdot(p-1)}\cdot scal_h\ .\]
    
We also want to discuss the case when $\alpha_-$ is a lightlike
nc-Killing form whose corresponding twistor is decomposable.
We use the following convention. If a decomposable $p$-form $\gamma$
on $\RR^{r,s}$ 
is the $\wedge$-product of lightlike $1$-forms only, i.e.,
the corresponding subspace to $\gamma$ in $\RR^{r,s}$ is totally 
lightlike 
then we call the decomposable $p$-form totally lightlike (isotropic)
as well.
There is a version of Lemma \ref{Lem11} for totally lightlike 
$p$-forms.
\begin{LM} \label{Lem12} Let $\alpha_-$ be a 
totally lightlike conformal Killing $p$-form 
on a space $(M,g)$ with the following two properties:
\begin{enumerate}
\item There is  $A\in \Gamma(TM)$ such that
$d\alpha_-=A^\flat\wedge\alpha_-$
and
\item
there is $B\in \Gamma(TM)$ such that $d^*\alpha_-=B\ecke\alpha_-$.
\end{enumerate}
Then it exists (locally) a rescaled metric $\tilde{g}$ in the
conformal class
$[g]$ such that the rescaled nc-Killing form $\tilde{\alpha}_-$ is
parallel. In particular, the holonomy
of the Levi-Civita connection to $\tilde{g}$ is
reducible with a fixed lightlike subspace. 
\end{LM}
{\bf Proof.} First, we show that we can assume   
$\alpha_-$ to be a closed form.  This is 
for the following reason. The differential 
form $d\alpha_-$
is decomposable and closed. Hence, by
Frobenius' there are (local) coordinates 
$(x_1,\ldots,x_n)$ such that $d\alpha_-=dx_1\wedge\ldots\wedge 
dx_{p+1}$.
Moreover, since $\alpha_-$ is decomposable we can choose these
coordinates such that $\alpha_-=f\cdot dx_1\wedge\ldots \wedge dx_p$,
where $f$ is a function in the $x_1,\ldots,x_{p+1}$. By rescaling the 
metric
with the function $f$ we find that $\tilde{\alpha}_-=f^{-1}\alpha_-$
is a closed nc-Killing form. 

Now let $\alpha_-=l_1\wedge\ldots\wedge l_p$ be a totally isotropic and 
closed conformal Killing form with $d^*\alpha_-=t\cdot l_1\wedge\ldots\wedge 
l_{p-1}$, 
where
the $l_i$'s are mutually orthogonal lightlike $1$-forms and $t$ is 
some 
function. 
Then we calculate in an arbitrary  point $m\in M$:
\begin{eqnarray*}
0&=&X(g(\bar{l}_1\ecke\ldots\ecke 
\bar{l}_{p-1}\ecke\alpha_-,\bar{l}_1\ecke\ldots\ecke 
\bar{l}_{p-1}\ecke\alpha_-))\\
&=&2\cdot g(\bar{l}_1\ecke\ldots\ecke 
\bar{l}_{p-1}\ecke\nabla_X\alpha_-,\bar{l}_1\ecke\ldots\ecke 
\bar{l}_{p-1}\ecke\alpha_-)\\
&=&\frac{2\cdot (-1)^p}{n-p+1}\cdot g(t
X^\flat,l_p)\qquad\mbox{for\ 
all}\ \ 
X\in T_pM\ ,\end{eqnarray*}
where we have chosen lightlike $1$-forms $\bar{l}_i$ with 
$\nabla\bar{l}_i(m)=0$ and 
\[g_m(l_i,\bar{l}_i)=1\qquad\mbox{and}\qquad g_m(l_i,\bar{l}_j)=0
\quad\mbox{for}\ \ i\neq j\ .\]
But this is only possible for all $X\in TM$ if $t\equiv 0$, i.e., 
$d^*\alpha_-=0$.
Henceforth, $\alpha_-$ is parallel and totally isotropic.
This proves the statements of the lemma.\hfill$\Box$\\

Using this lemma and the normal form description for decomposable twistors 
with respect to some $(s_-,s_+)$ leads us to the next result,
which is for lightlike nc-Killing forms.
\begin{LM} Let $\alpha_-$ be a totally isotropic 
nc-Killing $p$-form on $(M,g)$ with decomposable twistor. Then
there is (at least locally) a metric $\tilde{g}$ in the conformal class 
such that
the rescaled form $\tilde{\alpha}_-$ is parallel.  
\end{LM}
We can say even more than stated in the lemma.
With respect to the metric $\tilde{g}$, where
$\tilde{\alpha}_-$ is totally isotropic and parallel the 
corresponding twistor takes
the form
\[\alpha=(s_-^\flat+a\cdot s_+^\flat)\wedge l_1\wedge\ldots\wedge l_p\]
for some constant $a$.
However, if $a\neq 0$ then 
$\beta=l_1\wedge\ldots\wedge l_p$ would be a twistor itself, 
since the totally lightlike subspace, which uniquely belongs
to $\alpha$, is parallel with respect to $\nabla^{NC}$. This is 
not 
possible, because the fact that $\beta$ is a twistor 
means $\beta_-$ is zero and $d\beta_-\neq 0$. For this 
reason, the constant $a$ must be zero (so that $\alpha$
is totally isotropic), which implies 
that the scalar
curvature of $\tilde{g}$ is zero. Furthermore, the 
twistor
equations (\ref{g2}) and (\ref{g3}) show that the Ricci tensor
of $\tilde{g}$ maps into the totally lightlike subspace of the tangent 
space 
that corresponds
to the nc-Killing form $\tilde{\alpha}_-$. Then the metric $\tilde{g}$ has 
reducible holonomy with an invariant lightlike subspace
(that is not dilated under the action). And
the holonomy is possibly undecomposable 
(that is the generic case). The derived results so far
sum up to the following proposition.

\begin{PR} \label{PR100}
Let $(M,c)$ be a simply connected conformal space and $Hol(\omega_{NC})$ 
the corresponding holonomy group of the normal conformal 
connection $\omega_{NC}$. The holonomy group $Hol(\omega_{NC})$
fixes a decomposable $(p+1)$-form on $\RR^{r+1,s+1}$ 
($p=1,\ldots,n-1$) if and only if 
one of the following cases occurs.
\begin{enumerate}
%\item
%There is a Ricci-flat metric $g\in c$ and $Hol(\omega_{NC})$
%fixes a lightlike $1$-twistor.
%\item
%There is an Einstein metric $g\in c$ with $scal_g>0$ and 
%$Hol(\omega_{NC})$ fixes a timelike $1$-twistor.
%\item
%There is an Einstein metric $g\in c$ with $scal_g<0$ and
%$Hol(\omega_{NC})$ fixes a spacelike $1$-twistor.
\item
There is a product $h^p\times l^q$ of Einstein metrics in $c$
with \[scal_l=-\frac{q(q-1)}{p(p-1)}\cdot scal_h\ . \]
If $scal_h\neq 0$ then $Hol(\omega_{NC})$ fixes a non-degenerate subspace 
and if $scal_h=0$ then $Hol(\omega_{NC})$
fixes a degenerate subspace of dimension $p+1$ with exactly one lightlike 
direction.  
\item
There is $g\in c$ with totally isotropic Ricci 
tensor and parallel totally isotropic form.
The group 
$Hol(\omega_{NC})$
fixes a totally isotropic subspace (without dilation) of dimension 
at least $2$.
\end{enumerate}
\end{PR}

In both cases of the proposition the holonomy representation
$\RR^{r+1,s+1}$ is reducible. A further possibility for a reducible
holonomy $Hol(\omega_{NC})$ is the case where a lightlike subspace
is invariant, but dilated under the action. In this case the
volume form to the invariant subspace is not fixed by the action
of the holonomy, and therefore, does not give rise to a nc-Killing form, 
i.e., Proposition \ref{PR100} does not apply to this situation.

%%%%%%%%%%%%%%%%%%%%%%%%%%%%%%%%%%%%%%%%%%%%%%%%%%%%%%%%
\section{A geometric description for spaces with %%%%%%%
nc-Killing forms}                                %%%%%%%%
%%%%%%%%%%%%%%%%%%%%%%%%%%%%%%%%%%%%%%%%%%%%%%%%%%%%%%%%%
\label{p8}
Here we want to apply the results developed
in the two preceding paragraphs 
for conformal Einstein spaces and decomposable normal twistors
to derive a certain geometric description for conformal 
structures (metrics) admitting nc-Killing forms in case that the
normal conformal holonomy is decomposable.
For this description the reducibility of the stabiliser 
$S_\alpha:=Stab(\alpha)$
of a twistor $\alpha$ to a given nc-Killing form $\alpha_-$
is used to decompose the underlying conformal geometry.
The corresponding `(weakly) irreducible' parts are the building blocks of 
solutions.
In Theorem \ref{TH1} we discussed already examples for those
`irreducible' normal conformal geometries.

\begin{LM} \label{LM33} Let $Hol(\omega_{NC},c)$ be the
normal conformal holonomy of a simply connected space $(M,c)$. 
We assume that
there is a product of Einstein
metrics $h^p\times l^q$ in $c$ with $scal_l=-\frac{q(q-1)}{p(p-1)}\cdot
scal_h\neq 0$. Then it is
\[Hol(\omega_{NC},c)=Hol(\omega_{NC},[h])\times
Hol(\omega_{NC},[l])\]
and the holonomy representation decomposes
$\RR^{r+1,s+1}$ in two non-degenerate subspaces $V_1\oplus V_2$.
\end{LM}

{\bf Proof.} Over $M$ we have the $\SO(r+1,s+1)$-bundle 
$\frak{M}(M)$ with connection $\omega_{NC}$.
Let $H$ and $L$ denote the spaces where $h$ and $l$ live.
We can pull back the bundles
$\frak{M}(H)$ and $\frak{M}(L)$ with their normal connections to $M$
to obtain a $\SO(r_1+1,s_1+1)\times\SO(r_2+1,s_2+1)$-bundle with connection
$\omega$. The structure group of 
the latter bundle 
sits in $\SO(r+2,s+2)$. We denote the corresponding extended 
principal fibre bundle with $B(H,L)$.
Since $h$ 
and $l$ are
Einstein, the extended connection 
$\omega$ 
can be reduced to a subbundle $T$ of $B(H,L)$ with structure 
group $\SO(r+1,s+1)$.

We define now a bundle embedding of $\frak{M}(M)$
in $B(H,L)$. This embedding induces 
an 
isomorphism 
of the connection $\omega_{NC}$  
and the reduction of $\omega$ to the image of the embedding,
which is just the bundle $T$.
The map can be given with respect to a local frame
$(s_-,s_+,s_1,\ldots,s_n)$,
which fits to the scaling $g=h\times l$ such that 
$a=(s_1,\ldots,s_p)$ spans $TH$ and $b=(s_{p+1},\ldots,s_{n})$
spans $TL$, in the following way. Let
\begin{eqnarray*}
s_-^\flat+\frac{-scal_g}{2(n-2p)(n-1)}\cdot s_+^\flat 
&\quad\mapsto\quad&
a_-^\flat+\frac{scal_h}{2p(p-1)}\cdot a_+^\flat\\
s_-^\flat+\frac{scal_g}{2(n-2p)(n-1)}\cdot s_+^\flat&\quad\mapsto\quad&
b_-^\flat+\frac{scal_l}{2(n-p)(n-p-1)}\cdot b_+^\flat\\[3mm]
s_i&\quad\mapsto\quad& a_i,\qquad i=1,\ldots,p\\
s_i&\quad\mapsto\quad& b_i,\qquad i=p+1,\ldots, n\ .
\end{eqnarray*}

These assignments define a map for the whole
frame 
$(s_-,s_+,s_1,\ldots,s_n)$ to $B(H,L)$ and after extension by
right multiplication with $\SO(r+1,s+1)$ on the principal fibre bundles
we obtain the desired bundle map.
It is a straightforward calculation to see that 
\begin{eqnarray*} 
\omega_{NC}\circ ds_c(X)&=&\ \ \ s_-^\flat\wedge 
X^\flat+\omega_{LC}\circ
ds_c(X)-s_{+}^\flat\wedge K_g(X)^\flat\\ 
&=&\ \ \ a_-^\flat\wedge
X^\flat+\omega_{LC}\circ
da_c(X)-a_{+}^\flat\wedge K_h(X)^\flat\\
&&+\ b_-^\flat\wedge
X^\flat+\omega_{LC}\circ
db_c(X)-b_{+}^\flat\wedge K_l(X)^\flat\\
&=&\ \ \ \omega\circ d(a_c\!+\!b_c)(X)=\omega\circ 
da_c(X)+\omega\circ db_c(X)
\end{eqnarray*}
and this shows that
$Hol(\omega_{NC},c)$ is the product of the normal conformal holonomies
on $[h]$ and $[l]$. The condition on the subrepresentations $V_1$
and $V_2$ 
to be non-degenerate
is clear from Proposition \ref{PR100}.\hfill$\Box$\\

The subspaces $V_1$ and $V_2$ in the lemma are naturally identified
with the tangent spaces of the cones over $h$ and $l$ at every point
via the mapping 
given in the proof. 
There is also a version of Lemma \ref{LM33} when
$g=h\times l$ 
in $c$ is a product of Ricci-flat metrics (i.e., $g$ itself
is Ricci-flat). In this case it still holds
$Hol(\omega_{NC},c)=Hol(\omega_{NC},[h])\times
Hol(\omega_{NC},[l])$.

We extend now Proposition \ref{PR100} to a generalised form to make a  
statement
for the case when the holonomy representation decomposes
into an arbitrary number of non-degenerate components:
\[ \RR^{r+1,s+1}\ =\ \oplus_i V_i\ ,\qquad i=1,\ldots, v\ .\]
Indeed, then it is 
possible
to decompose the conformal space into further parts
and their normal conformal holonomies as well. However,
by applying this method one must pay attention to the fact that all the 
scalings
of the conformal structure when the metric on the base manifold becomes
a product in the conformal class will be different ones, in general. 
Hence, the conformal structure
is not just given by a simple product of several metrics. 
%We will come back to this point below in Theorem \ref{TH2}.
%Moreover, there is also the case when the decomposition
%$\RR^{r+1,s+1}=\oplus_i V_i$
%contains subrepresentations, which are degenerate with
%more than one lightlike direction. This implies
%that with respect to the right scaling there occur `components'
%in the conformal geometry, which
%are given by the conformal class of a weakly irreducible 
%pseudo-Riemannian metric. We will come back to these points below
%in Theorem \ref{TH2}. 
\begin{PR} \label{PR44} Let $(M,c)$ be a conformal space and
let $\RR^{r+1,s+1}\ =\ \oplus_i V_i$ 
be a decomposition of the representation of $Hol(\omega_{NC})$ into
(weakly) irreducible submodules 
$V_i$ with $dim V_i>1$ for all $i$. 
Then there exist (locally) metrics $g_i$, which are
Einstein with $scal_{g_i}\neq 0$ of signature $(r_i,s_i)$ and 
$Hol(\omega_{NC},[g_i])$
acts (weakly) irreducible on a subspace of codimension $1$ in 
$\RR^{r_i+1,s_i+1}$. Moreover,
there are functions $\phi_i$ such that 
\[c=\big[\sum_ie^{\phi_i}\cdot g_i\big]\ .\]
The holonomy group $Hol(\omega_{NC})$ is equal to the product
\[Hol(\omega_{NC},c)= \Pi_i\ Hol(\omega_{NC},[g_i])\ .\]
\end{PR}  

The statement of Proposition \ref{PR44} also implies that the cone
$\hat{g}_i$ of every factor $g_i$ in the decomposition has (weakly)
irreducible holonomy with respect to the Levi-Civita connection
(since $V_i$ can be identified with the tangent spaces of the
cone $\hat{g}_i$ at every point).

Now we can consider the situation when
a nc-Killing form $\alpha_-$ on a space $(M,c)$ exists such that the
stabiliser $S_\alpha$ of the
corresponding twistor $\alpha$ in $\SO(r+1,s+1)$
acts decomposable on $\RR^{r+1,s+1}$. In general, we have the
following decomposition.
\begin{LM} \label{LM44} Let $\alpha$ be a $p$-form on $\RR^{r+1,s+1}$
with stabiliser $S_{\alpha}\subset\SO(r+1,s+1)$ and let
\[\oplus_i V^\alpha_i=\RR^{r+1,s+1},\qquad i=1,\ldots,\nu \]
be the decomposition of the $S_\alpha$-module $\RR^{r+1,s+1}$
into (weakly) irreducible components (without dilation). Then there are 
unique
$p$-forms $\alpha_i\in\Lambda^p(V^\alpha_i)$ for all $i$, which
are stable under $S_\alpha$ such that
$\alpha=\sum_i\alpha_i$.
\end{LM}

Let $\alpha_i\neq 0$ for some $i$ as in Lemma \ref{LM44} 
with the property that $V_i^\alpha$ is (weakly) irreducible.
Then we obtain 
a factor 
$h$ in an appropriate scaled  metric $g=h\times l\in c$, and the conformal
structure $[h]$ admits a nc-Killing form $\alpha_{i-}$, which comes from
the twistor $\alpha_i$ on $V_i^\alpha$. Thereby, we identify
$V_i^\alpha$ with the tangent space of the cone over
$h$ at every point. The pullback 
$\alpha_{i-}$ to the total space $(M,g)$ is a nc-Killing form
itself and it holds
\[Y^\flat\wedge\alpha_{i\mp}=-\nabla_Y\alpha_{i-}+Y\ecke\alpha_{i0}=0\qquad
\mbox{for\ all}\quad Y\ \bot\ TH\ .\]
That implies $\alpha_{i\mp}=0$ and hence, $\alpha_{i-}$ is coclosed for 
$g$. It is also coclosed on $(H,h)$. We obtain the geometric
characterisation in case that the normal conformal holonomy is decomposable
(but the conformal structure not Einstein).

%Finally, we have to say something on the case that there is
%an irreducible component in $\RR^{r+1,s+1}$, which is 
%degenerate and contains more than one lightlike direction. 
%In this case there is a decomposition $\RR^{r+1,s+1}=\oplus_i V_i$
%with at least one $i$ such that $V_i$ is not irreducible,
%but has an invariant degenerate subspace with more than one 
%lightlike direction. Moreover, we assume that
%this subspace is not dilated, i.e., there is really some 
%lightlike decomposable $p$-form, which is invariant. The conclusion 
%is that there 
%is a scaling
%$g\in c$ such that $g=h$ or $g=h \times l$ with $h$ a weakly irreducible
%metric with totally isotropic Ricci-tensor.
%Taking together all the parts of our discussion here we arrive
%at the following description for conformal spaces 
%with nc-Killing forms.
\begin{THEO} \label{TH2} Let $\alpha_-$ be a nc-Killing form on $(M,c)$
and $\oplus_i V^\alpha_i$ ($dim V_i^\alpha>1$) a decomposition of 
$\RR^{r+1,s+1}$
with respect to the stabiliser $S_\alpha$ into (weakly) irreducible
components and let  $c=\big[\sum_i e^{\phi_i}\cdot 
g_i\big]$ 
be 
a 
corresponding representation
of the conformal class (as in Proposition \ref{PR44}). It is
$\alpha_-=\sum_i\alpha_{i-}$ and
if the component $\alpha_i$ on $V_i^\alpha$ to the twistor $\alpha$
is non-trivial then $\alpha_{i-}$ is a 
special Killing form for $g_i$. 
\end{THEO}

We remember to the fact that if $dim V_i^\alpha=1$ for some $i$
then the conformal structure is Einstein and we find a
characterisation of solutions via the ambient metric (cf.
Proposition \ref{PR50}, Theorem \ref{TH1}). In the next, paragraph we 
will 
come 
across the case
when $Hol(\omega_{NC})$ acts irreducible on $\RR^{r+1,s+1}$.
We also obtain the following conclusion concerning coclosed
nc-Killing forms.

\begin{CO} Let $\alpha_-$ be a nc-Killing form without zeros
on $(M,c)$ such that the stabiliser $S_\alpha$ decomposes $\RR^{r+1,s+1}$
to $V_1\oplus V_2$  with $dimV_{1,2}>1$. Then $\alpha_-$ is coclosed
with respect to the corresponding scaling $h\times l\in c$.
\end{CO}

\begin{table}[h!]
  \centering
  \setlength{\extrarowheight}{3pt}
  \renewcommand{\arraystretch}{1.3}
  \begin{tabular}{|>{$}c<{$}|>{$}c<{$}|>{$}c<{$}|}
    \hline
\RR^{r+1,s+1}=\oplus_i V_i& \mbox{scaling\ of\ } $c$ & \mbox{nc-Killing\ 
form\ } \alpha_-\\
    \hline\hline
 \begin{array}{c}V_i\ \mbox{non-degenerate},\\[-2.8mm] dimV_i>1\quad \forall\ 
i
\end{array} & \begin{array}{c}\sum_ie^{\phi_i}g_i,\\[-2.8mm]
g_i\ 
\mbox{Einstein}\end{array} & \begin{array}{c}
\alpha_-=\sum_i\alpha_{i-},\\[-2.8mm] \alpha_{-}\ \sim\ 
\mbox{coclosed}\end{array}\\
\begin{array}{c}V_i\ \mbox{non-degenerate},\\[-2.8mm] dimV_1=1\end{array}& 
\begin{array}{c} g\in c\ 
\mbox{Einstein},\\[-2.8mm] scal_g\neq 0\end{array} & 
\begin{array}{c}\alpha_-=\beta_-+\gamma_-,\\[-2.8mm] d^*\beta_-=d\gamma_-=0
\end{array}\\
\begin{array}{c} \mbox{undecomposable},\\[-2.8mm] \mbox{fixed\ 
1-form}\end{array} & 
g\in c\ \mbox{Ricci-flat} & 
\alpha_-\ \mbox{or}\ d\alpha_-\ \mbox{parallel}\\
\begin{array}{c} \mbox{undecomposable},\\[-2.8mm] \mbox{fixed\ p-form},\ 
1<p<n+1\end{array}  & \mbox{non-Einstein} &
\alpha_-\ \mbox{or}\ d\alpha_-\ \mbox{parallel}\\
\mbox{irreducible}& \mbox{non-Einstein} &\begin{array}{c} \alpha_-\ 
\mbox{nc-Killing}\\[-3mm] d^*\alpha_-=?\end{array}\\
\begin{array}{c}\mbox{undecomposable},\\[-2.8mm] \mbox{(fixed subspace\ with 
dilation)}\end{array}&
 \mbox{non-Einstein} & \mbox{non}
\\[3pt]
    \hline
  \end{tabular}
  \vspace{8pt}
  \caption{\label{Tabo} Normal conformal holonomy actions,
natural scalings of the underlying conformal structure and properties
of the occurring nc-Killing forms thereof.}
\end{table}

We can also understand now that if a nc-Killing form $\alpha_-$
without zeros on a non-conformal Einstein space
can not be scaled (locally) to a coclosed form
then the stabiliser $Stab(\alpha)$ of the corresponding twistor
acts irreducible or undecomposable on $\RR^{r+1,s+1}$.
In case that the space is conformally Einstein ($scal\neq 0$)
we mentioned already that every nc-Killing form is the sum
of a closed and a coclosed form in the Einstein scaling.
In Table \ref{Tabo} we give an overview of the situation that we
tried to explain here.

%%%%%%%%%%%%%%%%%%%%%%%%%%%%%%%%%%%%%%%%%%%%%%%%%%%%%%%%%%%%%%
\section{Solutions in dimension 4}%%%
%%%%%%%%%%%%%%%%%%%%%%%%%%%%%%%%%%%%%%%%%%%%%%%%%%%%%%%%%%%%%%
\label{p9}
After the general discussion so far, we want 
to study in this paragraph solutions 
of the normal twistor 
equations for differential
forms on $4$-dimensional Riemannian and Lorentzian manifolds
and their corresponding possible normal conformal holonomy groups.

\subsection{The Riemannian case}
Let $(M^4,g)$ be an oriented  Riemannian $4$-space.
We discuss nc-Killing forms 
according to their
degree. 

First, 
if there is a normal conformal function
$f_-=\alpha_-$ without zero then 
$M^4$ is a conformal Einstein space.
The $1$-form \[o=s_-^\flat-\frac{scal}{2(n-1)n}s_+^\flat\ \in\
\Omega^{1}_{\frak{M}}(M)\] is
parallel. That implicates that for $scal>0$ the holonomy
$Hol(\omega_{NC})$ 
of the normal conformal connection 
$\omega_{NC}$
is reduced at least to the subgroup $\SO(5)$ of the 
M{\"o}bius group $SO(1,5)$. For negative scalar
curvature ($scal<0$) the holonomy is reduced to 
a subgroup of $\SO(1,4)$. For Ricci-flat spaces the 
holonomy $Hol(\omega_{NC})$ is contained in the stabiliser
of a single lightlike vector in $\RR^{1,5}$. 

Next we take a look at the nc-Killing $1$-forms.
So let $\alpha_-$ be a nc-Killing $1$-form on $(M^4,g)$
with
$V_-=\alpha^\sharp_-$ its dual conformal vector field and
$f=g(V_-,V_-)$ the  square length. 
We assume for the moment that $f$ has no zero.
Then the conformally changed metric
$\tilde{g}=f^{-1}\cdot g$ has $V_-$ as a Killing vector field,
i.e.,
$L_{V_-}\tilde{g}=0$. Moreover, $V_-$ is divergence-free
and  the dual
$\tilde{\alpha}_-=\frac{1}{\sqrt{f}}\alpha_-$ 
is coclosed with respect to $\tilde{g}$. Now there are  
two
possibilities. 
Either $\tilde{\alpha}_-$ is parallel or it 
is not closed. In the first case 
one can see from the twistor equations that 
$\tilde{g}$
is Ricci-flat
with a parallel vector, i.e., $(M^4,g)$ is (conformally) flat.
 
In the latter case
we can find locally 
an orthonormal frame $\tilde{s}$ such that  
\[\tilde{\alpha}_-=\tilde{s}_1^\flat\qquad\mbox{and}\qquad
d\tilde{\alpha}_-=h\cdot \tilde{s}_2^\flat\wedge \tilde{s}_3^\flat\]
for some 
local function $h$. The integrability conditions 
(\ref{G17}) and (\ref{G18})  
then say that
\[ \tilde{s}_1\ecke W_{\tilde{g}}=0\qquad\mbox{and}\qquad 
\tilde{s}_1\ecke
(C(X,Y)\wedge\tilde{\alpha}_-)
=0\] 
for 
all $X,Y\in TM$. 
Moreover, it is
\[C_{\tilde{g}}(X,Y)\ecke\tilde{\alpha}_-=0\]
by (\ref{G19}), which together implies  
$C_{\tilde{g}}\equiv 0$. Therefore, 
we have $W_{\tilde{g}}(X,Y)\circ
d\tilde{\alpha}_-=0$ 
for all $X$ and $Y$, which means that 
$W_{\tilde{g}}(\tilde{s}_2^\flat\wedge
\tilde{s}_3^\flat)=k\cdot \tilde{s}_2^\flat\wedge \tilde{s}_3^\flat$ for 
some 
function 
$k$
and $W_{\tilde{g}}(\tilde{s}_i^\flat\wedge \tilde{s}_j^\flat)=0$ in all 
other 
cases. But
since $trW_{\tilde{g}}=0$,
the function 
$k$ must be zero and therefore, it is
$W_{\tilde{g}}\equiv 0$, i.e., $M^4$ is
conformally flat. Altogether, we can conclude that any $4$-space
admitting a nc-Killing $1$-form with or without zeros
is conformally flat.

Finally, we have to consider the case when 
$\alpha_-$ is a nc-Killing $2$-form. 
We assume without loss of generality that $\alpha_-$ is selfdual.
Moreover, we have seen already that
if $\alpha_-$ has no zeros, we can also
assume
that $\tilde{\alpha}_-$ is parallel with respect 
to some conformally changed metric
$\tilde{g}$ (cf. \cite{Sem01} and  paragraph \ref{p7}). Therefore, 
$\tilde{g}$
is a K{\"a}hler metric. There are two cases. First,
the metric $\tilde{g}$ is Einstein,
which is only possible if $\tilde{g}$ 
is Ricci-flat, i.e., $\tilde{g}$ has holonomy $\SU(2)$
or is flat. In the former case the holonomy of
$\omega_{NC}$ is the stabiliser group 
$Stab(e_-^\flat\wedge\omega_o)$,
where $\omega_o$ denotes the standard K{\"a}hler form
on $\RR^4$. 

In case that $\tilde{g}$ is not Einstein
it is locally up to
a constant scaling factor a product of the sphere 
$S^2$ with the hyperbolic space form $H^2$. 
The volume forms of the factors of this product
are the nc-Killing $2$-forms. The stabiliser in $\SO(1,5)$ 
of the corresponding twistors is $\SO(3)\times\SO(1,2)$. 
However, the product $S^2\times H^2$ is already conformally 
flat, i.e., the holonomy $Hol(\omega_{NC})$ is trivial.

\begin{THEO} \label{TH3}
Let $(M^4,[g])$ be a Riemannian $4$-space with
conformal structure $[g]$ and let $\alpha_-$ be 
a nc-Killing form without zero of arbitrary degree then 
at least one of the following cases occurs
(up to a conformal scaling factor)
\begin{enumerate}
\item
$deg(\alpha_-)=0$ and $M^4$ is Einstein,
\item
$deg(\alpha_-)=2$ and $M^4$ is Ricci-flat and K{\"a}hler or
%\item
%$deg(\alpha_-)=2$ and $M^4$ is locally $S^2\times H^2$ or
\item
$M^4$ is flat.
\end{enumerate}
In particular, 
the holonomy group of $\omega_{NC}$ for a simply connected  space $M^4$
is contained in one of the following subgroups of the M{\"o}bius group
$\SO(1,5)$:
\[\SO(5),\quad\SO(1,4),\quad Stab(e_-^\flat),\quad
Stab(e_-^\flat\wedge\omega_o)
\quad\mbox{or}\quad
\{e\}\ .\]
\end{THEO}
Obviously, irreducible representations
of subgroups of the M{\"o}bius group 
do not occur in the Riemannian case. 
It is a matter of further investigation, which proper subgroups 
of the stated groups in Theorem \ref{TH3} can really occur as holonomy 
groups.
For example, it is possible that there is a conformal geometry
such that $Hol(\omega_{NC})\subset \SO(1,4)$ 
acts weakly irreducible on $\RR^{1,4}$ with dilation. It is also 
possible
that $Hol(\omega_{NC})$ acts weakly irreducible on $\RR^{1,5}$
with dilation, although
in this case no nc-Killing form exists.
At least, we can say if $Hol(\omega_{NC})$ has an invariant subspace 
without dilation then the geometry is Einstein or conformally flat.
Table \ref{Tab11} gives an overview of possible holonomy
groups. 

\begin{table}[h!]
  \centering
  \setlength{\extrarowheight}{3pt}
  \renewcommand{\arraystretch}{1.3}
  \begin{tabular}{|>{$}c<{$}|>{$}c<{$}|>{$}c<{$}|}
    \hline
Hol(\omega_{NC})\ \mbox{sitting\ in} & \mbox{(local conformal)\ geometry}\ 
[g]& 
\mbox{nc-Killing\ form}\\
    \hline\hline
    \SO(5) & \mbox{Einstein},\ scal>0 & \mbox{one\ function}\\
   \SO(1,4) & \mbox{Einstein},\ scal<0 & \mbox{one\ function}\\
  Stab(e_-^\flat) & \mbox{Ricci-flat} & \mbox{one\ function}\\
    Stab(e_-^\flat\wedge\omega_o) & 
\mbox{Ricci-flat,\ K{\"a}hler} & 
\mbox{K{\"a}hler form}\\
%\SO(3)\times \SO(1,2) & S^2\times H^2 & 2\times\ \ 2\mbox{-form} \\
\{e\} & \mbox{conformally\ flat} & \mbox{maximal}\\
Stab(\RR\cdot e_-^\flat)
& ? & \mbox{non}
\\
SO(1,5) & \mbox{generic\ case} & \mbox{non} \\[3pt]
    \hline
  \end{tabular}
  \vspace{8pt}
  \caption{\label{Tab11} Partial holonomy list for the normal
conformal connection of Riemannian spaces in dimension $4$.} 
\end{table}

As we mentioned before it is well known
that there exits a 
selfdual nc-Killing $2$-form
with zero on the conformal completion of the Eguchi-Hanson metric,
which is not any longer conformally equivalent to an Einstein metric
(cf. \cite{KR96}). Nevertheless, the holonomy $Hol(\omega_{NC})$
in this case is equal to $Stab(e_-^\flat\wedge\omega_0)$.

\subsection{ The Lorentzian case}
We turn now to the case of a Lorentzian 
$4$-manifold $(M^{4,1},g)$. We
choose the 
signature $(-+++)$. If there is a nc-Killing function 
without zero
the space $M^{4,1}$ is conformally Einstein. 
In case that there is a 
timelike or spacelike nc-Killing $1$-form $\alpha_-$
one can easily show (as in the Riemannian case)
that $M^{4,1}$ is already conformally flat.

Here in the Lorentzian case we must also take into 
consideration the case
when the length $f:=g(\alpha_-,\alpha_-)$
of a nc-Killing $1$-form $\alpha_-$ vanishes identically,
i.e., the dual conformal vector field $V_-$
is everywhere null. We consider such a field locally and without
zeros. Then we can  
assume that 
$g$ is scaled such that $V_-$ is a Killing 
vector field on $M^{4,1}$. There
are two possible
cases. Either \[\alpha_-\wedge d\alpha_-\neq 0\] 
or 
this twist $3$-form
vanishes identically.
If it vanishes then $V_-$ is parallel in the scaling $g$. 
In particular, in dimension $4$
that means $g$ is a pp-wave with vanishing
scalar curvature. The corresponding twistor takes the form
$e_-^\flat\wedge\alpha_-$ determining  
a totally lightlike plane in $\RR^{2,4}$.
 
In the other case the twist $\alpha_-\wedge d\alpha_-$ does not 
vanish. This case is also well known. The underlying
metric $g$ is   
a so-called Fefferman metric (cf. \cite{Fef76}, \cite{Gra87}).
Equivalently, a Fefferman metric is determined by
the existence of a lightlike nc-Killing vector field $V_-$ with 
the property
\[Ric(V_-,V_-)=const>0\ .\]

\begin{table}[h!]
  \centering
  \setlength{\extrarowheight}{3pt}
  \renewcommand{\arraystretch}{1.3}
  \begin{tabular}{|>{$}c<{$}|>{$}c<{$}|>{$}c<{$}|}
    \hline
Hol(\omega_{NC})\ \mbox{sitting\ in} & \mbox{local\ conformal\ geometry}\
[g]&
\mbox{nc-Killing\ form}\\
    \hline\hline
    \SO(1,4) & \mbox{Einstein},\ scal>0 & \mbox{one\ function}\\
   \SO(2,3) & \mbox{Einstein},\ scal<0 & \mbox{one\ function}\\
Stab(e_-^\flat) &  \mbox{Ricci-flat} & \mbox{one\ function}\\
\begin{array}{c}Stab(e_-^\flat\wedge l),\\[-3mm]
l\ \mbox{null},\ \ l\bot e_-^\flat\end{array} &
\mbox{pp-wave} & \begin{array}{c}1\mbox{-form\ without\ twist},\\[-3mm]
Ric(V_-,V_-)=0\end{array} \\
\SU(1,2) & \mbox{Fefferman\ spaces} & 
\begin{array}{c}1\mbox{-form\ with\ twist},\\[-3mm] Ric(V_-,V_-)
> 0\end{array}\\
%\SO(3)\times \SO(2,1) & H^2\times S^{1,1} &
%2\times\ 2\mbox{-form} \\
%\SO(1,2)\times \SO(1,2) & S^2\times H^{1,1} &
%2\times\ 2\mbox{-form} \\
\{e\} & \mbox{conformally\ flat} & \mbox{maximal}\\
Stab(3\mbox{-form}) & ? & 2\mbox{-form} \\
\begin{array}{c}\mbox{undecomposable}\\[-3mm]
\mbox{with dilation}\end{array} &
? & \mbox{non}
\\
\SO(2,4)& \mbox{generic\ case} & \mbox{non}\\[3pt]
    \hline
  \end{tabular}
  \vspace{8pt}
  \caption{\label{Tab22} Possible holonomy groups for the normal
conformal connection of a Lorentzian space in dimension $4$.}
\end{table}

Fefferman spaces are locally constructed
as 1-dimensional fibre bundles over strict
pseudoconvex CR-manifolds $(N^3,W,J,\gamma)$,
where the direction of the fibre is that of $V_-$, i.e.,
lightlike, the $2$-form $d\alpha_-$ is
related to the complex structure $J:W\to W$ on
the CR-manifold $N$ and the $1$-form $\alpha_+$ is related
to $\gamma$, which is dual to the Reeb vector field of
the pseudoconvex CR-structure. 
The
normal twistor $\alpha$ in $\Omega_{\frak{M}}^2(M)$
belonging to $\alpha_-$
takes a very natural form. It corresponds to
the standard (pseudo)-K{\"a}hler form
$\omega_o$ on $\RR^{2,4}$. This implies that the holonomy of
$\omega_{NC}$ is included in $\SU(1,2)$
for every  Fefferman space. The group
$\SU(1,2)$ acts irreducible on $\RR^{2,4}$.
This characterisation   
by the holonomy group can 
be seen as a conformally invariant definition
for Fefferman spaces. It directly implies that if
a Fefferman space has the 
full 
holonomy group $\SU(1,2)$ it is not conformally Einstein.

Finally, we shortly mention some statements about 
the situation  for nc-Killing $2$-forms
in the Lorentzian setting.
The spaces
$S^2\times H^{1,1}$ and $H^{2}\times S^{1,1}$ have
parallel and decomposable nc-Killing $2$-forms, which
are the volume forms of the factors in the product. The 
stabilisers 
of these are $\SO(1,2)\times\SO(1,2)$ in the first
and $\SO(3)\times\SO(2,1)$ in the second case.
However, these product spaces are again conformally flat, i.e.,
if the holonomy $Hol(\omega_{NC})$ is contained in those 
stabilisers then it is trivial. 
But there may appear other $3$-twistors on a Lorentzian 
$4$-space. To investigate this situation, it could
be useful to know that there  
exists a complete orbit type classification
of the space $\Lambda^3\RR^6$ of $3$-forms
under the action of $\GL(6)$ (cf. \cite{Rei07}, \cite{Hit00}).
From these orbit types one can try to calculate the
orbit types with respect to the action of the
structure group $\SO(2,4)$.
And then one could investigate whether further
stabiliser groups of
those $3$-form orbits appear that are the holonomy groups to some special  
normal conformal geometries. In Table \ref{Tab22}
we list the cases for nc-Killing forms 
and holonomies of $\omega_{NC}$ for a Lorentzian
$4$-manifold that we mentioned here in our discussion. 

%%%%%%%%%%%%%%%%%%%%%%%%%%%%%%%%%%%%%%%%%%%%%%%%%%%%%%%%
\section{Application to conformal Killing spinors}  %%%%
%%%%%%%%%%%%%%%%%%%%%%%%%%%%%%%%%%%%%%%%%%%%%%%%%%%%%%%%
\label{p10}
The discussion of the normal twistor equations
for differential forms so far shows that 
there are methods to describe their solutions and many of their
underlying geometric structures are well known objects
and do occur in the mathematical literature as subject
to substantial work. We want to use our acquired
results to study a topic in conformal differential
geometry, which itself was subject to various investigations
during the last 30 years, 
namely  the twistor equation
for spinor fields.   
There is  a systematic investigation
of this twistor equation in Riemannian spin geometry
(cf. \cite{Lic88}, \cite{BFGK91}).
However, the origin of the twistor equation for spinors
was in the theory of General Relativity (cf. \cite{Pen67}, 
\cite{PR86}) and nowadays, 
there 
has been done considerable work 
for twistor spinors in the Lorentzian setting, too (cf. \cite{Bau00}, 
\cite{BL03}). 
The remarkable point of the twistor equation
for spinors, and this can explain why this equation
plays an important role, is the fact that all its 
solutions are automatically
normal (in the sense that we used here already). 
We start with recalling the very basic facts about spinors and
their twistor equation. 
   
Let $(M^n,g)$ be a semi-Riemannian spin manifold of dimension $n \ge 3$. We 
denote by $S$
 the spinor bundle and by $\mu :T^*M \otimes S \to S$ the
Clifford multiplication. The 1-forms with values in the spinor bundle
 decompose into two subbundles
\[
T^* M \otimes S = V \oplus Tw,
\]
where $V$, being the orthogonal complement to the `twistor bundle' $Tw := 
Ker \mu\;$,
is isomorphic to $S$.
We obtain two differential operators of first order by composing the spinor
derivative $\nabla^S$ with the orthogonal projections onto each of these
subbundles,
the Dirac operator $D$
\[
D: \Gamma (S) \stackrel{\nabla^S}{\longrightarrow} \Gamma (T^* M \otimes S)=
\Gamma (S \oplus Tw) \stackrel{pr_S}{\longrightarrow}\Gamma (S)
\]
and the twistor operator $P$
\[
P: \Gamma (S) \stackrel{\nabla^S}{\longrightarrow} \Gamma (T^*M \otimes S) =
\Gamma (S \oplus Tw) \stackrel{pr_{Tw}}{\longrightarrow} \Gamma (Tw).
\]
Both operators are conformally covariant. More exactly, if
$\tilde{g} = e^{-2 \phi } g$ is a conformal change of the metric,
the Dirac and the twistor operator satisfy
\begin{eqnarray*}
D_{\tilde{g}} &=& e^{ \frac{n+1}{2} \phi} D_g e^{-\frac{n-1}{2} \phi}\\
P_{\tilde{g}} &=& e^{ \frac{\phi}{2}} P_g e^{ \frac{\phi}{2}}.
\end{eqnarray*}

A spinor field is called conformal
Killing spinor if it lies in the kernel of the twistor operator
$P$. Alternatively, a spinor field $\varphi_- \in \Gamma(S)$ is a conformal 
Killing 
spinor if
and only if
\[
\nabla^S_X \varphi_- + \frac{1}{n} X \cdot D\,\varphi_- = 0  \qquad  
\mbox{for 
all vector fields}
\;\; X.
\]
This is the twistor equation for spinors.
A conformal Killing spinor $\varphi_-$ with respect to $g$ rescales
by \[\tilde{\varphi}_-:=e^{-1/2\phi}\cdot\varphi_-\] to a conformal 
Killing 
spinor 
with 
respect
to the conformally changed metric $\tilde{g}:=e^{-2\phi}\cdot g$. 
Thereby, we use the canonical identification of the spinor bundles
over $(M,g)$ and $(M,\tilde{g})$. We  
say that $\tilde{\varphi}_-$ is conformally equivalent to $\varphi_-$.

Obviously, each parallel spinor ($\nabla^S\varphi_-=0$) is a twistor
spinor. Another special class of twistor spinors are the Killing
spinors $\varphi_-$, which satisfy the equation \[\nabla_X^S \varphi_- = 
\lambda
X\cdot\varphi_-\] for all $X\in TM$ and some fixed 
$\,\lambda\in\CC 
\setminus \{0\}$ 
(this can be seen as corresponding notion to Killing 
differential forms).
The
dimension of the space of twistor spinors is a conformal invariant
and bounded by
\[
\dim \ker P \le 2 \,\mathrm{rank} S = 2^{[\frac{n}{2}]+1} =: d_n.
\]
If $\dim \ker P =d_n$, then $(M^n,g)$ is conformally flat.
Conversely, if $(M^n,g)$ is simply connected and conformally
flat, then $\dim \ker P= d_n$.

Now, we focus our attention on the case of Lorentzian signature $(-+
\ldots +)$. 
Let $(M^n,g)$ be an oriented and time-oriented Lorentzian spin
manifold. On the spinor bundle $S$ there exists an indefinite non-degenerate
inner product $\langle \cdot ,\cdot \rangle$ such that
\begin{eqnarray*}
\langle X \cdot \varphi_-  , \psi_- \rangle & =& \langle \varphi_- , X \cdot 
\psi_-
\rangle \quad \quad \mbox{and}\\
X(\langle \varphi_- , \psi_- \rangle )&= &\langle \nabla^S_X \varphi_- , 
\psi_-
\rangle + \langle \varphi_- , \nabla^S_X \psi_- \rangle,
\end{eqnarray*}
for all vector fields $X$ and all spinor fields $\varphi_-, \psi_-$
(cf. \cite{Bau81}). Each
spinor field $\varphi_- \in \Gamma (S)$ defines a vector field 
$V_{\varphi-}$
on $M$, the so-called Dirac current, by the relation
\begin{equation*}
g (V_{\varphi-} , X) := - \langle X \cdot \varphi_- , \varphi_- \rangle 
\end{equation*}
for all $X\in TM$. The vector $V_{\varphi-}$ is causal and future-directed. 
The zero sets of
$\varphi_-$ and $V_{\varphi-}$ coincide, i.e., for a non-trivial
spinor the associated field is non-trivial. (This is a very 
useful fact specific for Lorentzian geometry.) 
If $\varphi_-$ is a 
twistor
spinor, then $V_{\varphi-}$ is a conformal Killing field. The dual
of $V_{\varphi-}$ is denoted by $\alpha_{\varphi-}$. 

We have the 
following known geometric structure result for Lorentzian spaces with 
conformal
Killing spinors. Thereby, we call a space
with a parallel lightlike vector field a Brinkmann space. The notion
of Fefferman spaces that appeared in the preceding section can be extended
to every even-dimensional manifold ($n>2$) (cf. 
\cite{Gra87}, \cite{Bau98}). We just say 
here in 
short, 
a Fefferman
space $(M^{2m},[g])$ is a Lorentzian space with a conformal structure such 
that the normal conformal holonomy group $Hol(\omega_{NC},[g])$
sits in $\SU(1,m)$. A Lorentzian Einstein-Sasaki structure on
an odd-dimensional manifold $(M^{2m+1},g)$ lifts
to a K{\"a}hler structure on its Ricci-flat cone, i.e., the holonomy
$Hol(\omega_{NC})$ sits in $\SU(1,m)$, which is itself a subgroup of 
$\SO(2,2m+1)$.

\begin{PR} (\cite{BL03}) 
\label{PR7} Let $(M^n,g)$ be a Lorentzian spin manifold 
and 
let 
$\varphi_-\in\Gamma(S)$ be a conformal Killing spinor without zeros and 
let $V_{\varphi-}$ be 
its associated vector field. 

a) In case that $V_{\varphi-}$ is 
lightlike
there are exactly two different cases: 
\begin{enumerate} 
\item The twist $3$-form
$t_{\varphi-}:=\alpha_{\varphi-}\wedge d\alpha_{\varphi-}$ vanishes 
identically
and $\varphi_-$ is (locally) conformally equivalent to a parallel spinor, 
whose
associated lightlike vector field defines a Brinkmann space.
\item The twist $t_{\varphi-}$ does not vanish and $(M^n,g)$ is a Fefferman 
space.
\end{enumerate}

b) In case that $g(V_{\varphi-},V_{\varphi-})=const<0$ and 
$\varphi_-$ is a Killing spinor (i.e., $D\varphi_-=-n\lambda\varphi_-$) 
the space $(M^n,g)$ is Einstein-Sasaki.

c) If the length function $\langle\varphi_-,\varphi_-\rangle$ has no zero
then 
the 
metric
$g$ can be rescaled to the Einstein metric 
$\tilde{g}:=\frac{1}{\langle\varphi_-,\varphi_-\rangle^2}\cdot g$. 
In particular, 
there exists at least one Killing spinor on $(M,\tilde{g})$. 
\end{PR}

The results in Proposition \ref{PR7} do not give a complete answer to the 
case when the associated field $V_{\varphi-}$ is timelike. It is
only mentioned the Einstein-Sasaki case. It is our 
interest
to handle the geometric description for all timelike cases of
$V_{\varphi-}$. 
Thereby, 
both 
cases that $(M^n,g)$ is conformally Einstein resp. is not 
conformally Einstein are possible and of interest. 
For this purpose we can apply now our results from the discussion
of normal conformal Killing forms that we have developed in the previous
paragraphs. That this approach is reasonable is justified by the fact that,
as we mentioned already, every conformal Killing spinor is normal.
We can see this in the following way (for arbitrary signature). 

Assume that $(M^n,g)$ is a spin manifold. Then there is a spin M{\"o}bius
bundle on $M$, which we denote by $\frak{M}_{Spin}(M)$. This is a 
principal 
fibre
bundle with structure group $\Spin(r+1,s+1)$ and is a double cover
of $\frak{M}(M)$ respecting the right multiplication on the fibres
and the natural 
homomorphism \[\lambda:\Spin(r+1,s+1)\to \SO(r+1,s+1)\ .\]  
The normal conformal connection form $\omega_{NC}$ on $\frak{M}(M)$
lifts to a unique connection form on the spin M{\"o}bius bundle 
$\frak{M}_{Spin}(M)$.
Moreover, this connection form induces 
a covariant derivative $\nabla^{NCS}$ on the tractor spinor bundle 
$S_\frak{M}$ defined 
by
\[S_{\frak{M}}:=\frak{M}_{spin}(M)\times_\rho\Delta_{r+1,s+1}\ ,\]
where $(\Delta_{r+1,s+1},\rho)$ is the spinor representation
in signature $(r+1,s+1)$. 

With respect to the metric $g$, this spinor
bundle splits into the sum of two usual spinor bundles over $(M^n,g)$,
$S_{\frak{M}}=S\oplus S$.
The equation $\nabla^{NCS}\varphi=0$
translates to
\[\begin{array}{rcl} \nabla^S_X\varphi_-+X\cdot\varphi_+&=& 0\\[3mm]
\nabla_X^S\varphi_+-\frac{1}{2} K(X)\cdot\varphi_- &=& 0\ ,\end{array} \]
where $\varphi_-$ and $\varphi_+$ are smooth spinor fields on $(M,g)$.
From the two twistor equations, it follows that 
$\varphi_+=\frac{1}{n}\cdot D\varphi_-$ and 
the normal twistor equations for spinors $\varphi_-\in\Gamma(S)$ 
take the following form:
\begin{eqnarray} \nabla^S_X\varphi_-+\frac{1}{n}X\cdot 
D\varphi_-&=& 
0 \label{Spin1}\\[3mm]
\nabla_X^SD\varphi_--\frac{n}{2} K(X)\cdot\varphi_- &=& 0\ .\label{Spin2}
\end{eqnarray} 
We recognise that the first of the two equations is just
the conformal Killing spinor equation as introduced before.

Until now everything works analogous to the normal twistor equations for
differential forms. However, it is easy to see that the 
second equation for spinors is implied by the first equation alone,
the conformal Killing spinor equation. This is in contrast
to the case of differential forms, where the additional normal
twistor equations are not implicated by the conformal Killing equation.
Here for spinors we  
calculate:
\[\begin{array}{rcl}
\nabla^S_{s_i}\nabla^S_{s_j}\varphi_-+
\frac{1}{n}s_j\cdot\nabla^S_{s_i}D\varphi_-&=&0\\[4mm]
\nabla^S_{s_j}\nabla^S_{s_i}\varphi_-+
\frac{1}{n}s_i\cdot\nabla_{s_j}D\varphi_-&=&0\ ,\end{array}\]
which results to
\[ R^S(s_j,s_i)\cdot\varphi_-=-\frac{1}{n} s_i\cdot\nabla_{s_j}D\varphi_-+
\frac{1}{n}s_j\cdot\nabla_{s_i}D\varphi_-\ .\]
Using that $\sum_k s_k\cdot R^S(s_i,s_k)\varphi_-=
-\frac{1}{2} Ric(s_i)\varphi_-$ and 
$D^2\varphi_-=\frac{n\cdot scal_g}{4\cdot (n-1)}\cdot\varphi_-$, it 
follows the second twistor equation (\ref{Spin2}) for spinors.
\begin{THEO} \label{TH4} Let $(M,c)$ be a conformal spin space of 
signature 
$(r,s)$. The sets of 
conformal 
Killing spinors $\varphi_-\in\Gamma(S)$ (i.e., 
$\nabla^S_X\varphi_-+\frac{1}{n}X\cdot D\varphi_-=0$)
and normal twistor spinors $\varphi\in\Gamma(S_{\frak{M}})$ 
(i.e., $\nabla^{NCS}\varphi=0$) are 
naturally identified by the mapping
\[\varphi_-\ \mapsto\ \ \varphi\ =\ (\ \varphi_-\ ,\ \frac{1}{n}
\cdot D\varphi_-)\ 
.\]
\end{THEO}

In general, a spinor field gives rise to differential forms
of degree $p$. In Lorentzian
geometry we have already introduced as special case the Dirac current. 
The general construction is 
as follows. Let $\varphi_-$ be a spinor field on $(M^n,g)$. 
The corresponding
$p$-forms $\alpha^p_{\varphi-}$ are defined by the relation
\[g(\alpha^p_{\varphi-},X^p):=
i^{p(p-1)+r+1}\langle X^p\cdot\varphi_-,\varphi_-\rangle\qquad
\mbox{for\ all\ } X^p\in\Lambda^p(M)\ .\]
The so defined $p$-forms are not non-trivial in general. 
This depends
on the given spinor.
The same method applies to attach tractor $(p+1)$-forms 
to
some tractor spinor $\varphi$. Moreover, the tractor forms
associated to a twistor spinor are parallel with respect
to $\nabla^{NC}$, i.e., they are 
normal twistor forms. The corresponding induced nc-Killing forms
are just the associated differential forms of the induced conformal
Killing spinor. 
We want to direct 
our attention
to the Lorentzian case again.
Our
discussion 
gives rise to the following corollary to Theorem \ref{TH4}.
\begin{CO} Let $\varphi_-$ be a conformal Killing spinor
on a Lorentzian spin manifold $(M^n,g)$. Then the 
Dirac current $V_{\varphi_-}$ is a non-trivial and causal
nc-Killing vector field. Its dual $\alpha_{\varphi-}$ is a nc-Killing
$1$-form. 
\end{CO}
This corollary is now our starting point for the geometric
description of Lorentzian spaces admitting conformal
Killing spinors, in particular, those with timelike Dirac current. 

We have learned in paragraph \ref{p8} that the twistor
to an nc-Killing form either has an irreducible acting stabiliser 
on a subspace of codimension $0$ or $1$ in $\RR^{r+1,s+1}$ 
or else, for example
in the non-degenerate case, there exists a product metric in 
the conformal class. We want to use this philosophy here and show
that the rank $rk(\alpha_{\varphi-})$ of the nc-Killing $1$-form 
$\alpha_{\varphi-}$ to 
a conformal Killing spinor $\varphi_-$
determines whether or not there is (in the non-degenerate case)
a product
in the conformal class with respect to which the nc-Killing $1$-form
$\alpha_{\varphi-}$ restricts to the factors. Thereby, the rank 
$rk(\alpha_-)$ of an arbitrary 
$1$-form $\alpha_-$ is defined to be the 
unique natural number 
such that 
\[\alpha_-\wedge(d\alpha_-)^{rk(\alpha_-)}\neq 0\qquad\mbox{and}\qquad
\alpha_-\wedge(d\alpha_-)^{rk(\alpha_-)+1}=0\ .\]    
We notice that for an nc-Killing $1$-form $\alpha_-$ and all 
natural numbers
$l$ it is
\[(\alpha^{l+1})_-=\alpha_-\wedge(d\alpha_-)^l\ ,\]
i.e., the form $\alpha_-\wedge(d\alpha_-)^l$ is itself an nc-Killing form
and it corresponds to the twistor $\alpha^{l+1}$. 

The next point that is important for our discussion
is the fact that there exists a complete normal form classification for
skew-adjoint endomorphisms on the pseudo-Euclidean space $\RR^{2,n}$ with 
signature $(2,n)$ (cf. \cite{Bou02}). These normal forms can be written
down explicitly and their stabiliser groups in $\SO(2,n)$
can be calculated. We are 
interested here in those skew-adjoint endomorphisms whose
stabiliser is maximal (in the sense that the stabiliser is not 
properly contained in the stabiliser of any other skew-adjoint operator).
The space of skew-adjoint endomorphisms is naturally identified
with the space of $2$-forms on $\RR^{2,n}$. We present in the
following list all normal forms for $2$-forms on $\RR^{2,n}$ 
with maximal stabiliser and with the additional condition
that the corresponding skew-adjoint operators map
all causal vectors again to causal vectors. The latter requirement
characterises those $2$-forms, which are associated to
a spinor in the module $\Delta_{2,n}$ for signature $(2,n)$
(cf. \cite{Lei03}). 
\begin{enumerate}
\item $l_1\wedge l_2$ with $l_1$ and $l_2$ spanning a totally
lightlike subspace in $\RR^{2,n*}$, 
\item
$l_1\wedge t_1$ with $l_1$ lightlike and $t_1$ an orthogonal
timelike $1$-form to $l_1$,
\item
$\omega_o$ the standard symplectic form on $\RR^{2,n}$,
\item
$\omega_o|_V$ the standard symplectic form on 
a non-degenerate subspace $V$ of $\RR^{2,n}$ with signature $(2,p-1)$, 
where $p<n$.
\end{enumerate}
Moreover, we can say that every $2$-form $\alpha$ determines a unique
normal form with maximal stabiliser, which contains $Stab(\alpha)$. 

We use these normal forms in the following way. Let us assume
that $(M^n,g)$ is a Lorentzian spin manifold with conformal Killing 
spinor
$\varphi_-$. The corresponding twistor spinor 
$\varphi\in\Gamma(S_{\frak{M}})$ is parallel
with respect to $\nabla^{NCS}$ and induces a non-trivial twistor
$2$-form $\alpha_\varphi$.
The holonomy of the connection $\omega_{NC}$ lies 
necessarily in the stabiliser belonging to the normal form of the
$2$-twistor $\alpha_\varphi$. This implies that
there is a twistor $2$-form $\alpha$ on $M$ whose corresponding
normal form is one of those in the list, since they have
maximal stabiliser. Let us discuss 
the different
cases.

First, we assume that
the normal form to $\alpha$ is $l_1\wedge l_2$. In this case
the corresponding nc-Killing $1$-form is lightlike and hypersurface
orthogonal ($\alpha^2=0$). We can conclude that the underlying
conformal structure is represented by a Brinkmann metric.
However, for the particular case $l_1\wedge l_2$, one can 
even show that 
$\alpha_\varphi=\alpha$
(cf. \cite{Lei03}).
This implies that the twistor spinor $\varphi_-$ is 
(locally) conformally
equivalent to a parallel spinor on the Brinkmann space, which
induces the lightlike parallel vector.
In the second case, we can conclude, since $\alpha^2=0$, that
the underlying conformal structure is that of a static spacetime
with parallel spinor, i.e., it is  
\[[g]=[-dt^2+h]\ ,\] where $h$ is a Riemannian
metric with parallel spinor. The third case of a 
symplectic form 
$\omega_o$ is the 
Fefferman
case when the rank of $\alpha_-$ is $(n/2)-1$. It remains to investigate
the forth case when the normal form 
is the standard symplectic form on a proper subspace
of $\RR^{2,n}$. This is exactly the unknown case that we aim 
to 
describe here.

So let us assume that $(M^n,g)$ is a simply connected Lorentzian 
spin manifold
with a conformal Killing spinor $\varphi_-$, whose
corresponding twistor spinor $\varphi$ induces a twistor 
$2$-form
$\alpha_\varphi$
that has a stabiliser in \[\U(1,\frac{p-1}{2})\times\SO(n-p+1)\ .\] Then 
the $1$-form $\alpha_{\varphi-}$ is timelike and 
there exists a  
$2$-twistor
$\alpha$ with a restricted symplectic form 
as normal 
form, which induces a nc-Killing form $\alpha_-$ on $M$ of rank 
\[rk(\alpha_-)=rk(\alpha_{\varphi-})=\frac{p-1}{2}\ .\]
In particular, the twistor $\alpha^{rk(\alpha_-)+1}$ is
decomposable of degree $p+1$. The corresponding
subspace of $\RR^{2,n}$ is non-degenerate and we can conclude
that there is a product metric $h\times l$ in $[g]$, where
$h$ and $l$ are Einstein metrics on simply connected spin manifolds 
$H^{p}$ 
and
$L^{n-p}$ with dimension $p$ 
resp. $n-p$. Moreover,
$\alpha_-$ is a coclosed nc-Killing $1$-form of maximal rank 
$\frac{p-1}{2}$ on
$(H,h)$, and this implies that $h$ is an Einstein-Sasaki metric.  
In fact, $(H,h)$ is a simply connected Einstein-Sasaki space admitting
a conformal Killing spinor $\psi_{H-}$, which induces the nc-Killing 
$1$-form $\alpha_-$. However, the existence of the conformal Killing 
spinor
$\varphi_-$ should also impose a condition on $(L,l)$, and therefore, it 
still remains to 
discuss the 
geometry of the Riemannian
Einstein spin manifold $(L,l)$.      

For this we observe the following. The representation
$\RR^{2,n}$ splits into the subspaces $V$ and $W$ under the projection of 
the
stabiliser group
$Stab(\varphi)$ to $\SO(2,n)$. On $V$ lives the symplectic 
form $\omega_o$ with signature $(2,p-1)$.
Let $\Delta_V$ and $\Delta_W$ denote the 
spinor modules
over $V$ resp. $W$. 
As representations spaces of $\Spin(2,p-1)\times\Spin(n-p+1)$, it is 
\[\begin{array}{l} \Delta_V\otimes\Delta_W=\Delta_{2,n}\qquad\mbox{for}
\ n\ \mbox{odd},\\[1.5mm]
\Delta_V\otimes\Delta_W=\Delta_{2,n}^{\pm}\qquad\mbox{for}\ n\ 
\mbox{even}, 
\end{array}\]
where $\Delta^{\pm}$ 
denotes the half spinor modules
in even dimensions. 
The splitting with respect to $Stab(\varphi)$ gives
rise to a decomposition of
the tractor spinor bundle as
\[S_\frak{M}=S_V\oplus S_W \ .\]

We observe now that if we choose the product metric
$h\times l$ on $M=H\times L$ in the  conformal class 
$[g]$ then 
we can naturally identify
$S_V$ and $S_W$ with the spinor bundles over the cones 
$(\hat{H},\hat{h})$
resp. $(\hat{L},\hat{l})$ restricted to the bases:
\[S_V|_H\cong S_{\hat{H}}|_H\qquad\mbox{and}\qquad
S_W|_L\cong S_{\hat{L}}|_L\ .\]
We also know that it holds
\begin{eqnarray}Hol(\omega_{NC})= Hol(\hat{h})\times Hol(\hat{l})\ . 
\label{For1}\end{eqnarray}

Moreover, we have the following general fact.
Let $\rho$ be a representation of a product group $G_1\times
G_2$ and $\rho_1,\rho_2$ be
representations of $G_1$ resp. $G_2$ such that
\[\rho\cong\rho_1\otimes\rho_2\ .\]
Then the representation $\rho$ has a fixed vector if and only if
$\rho_1$ and $\rho_2$
both have fixed vectors. 
In our situation that means the
stabiliser
$Stab(\varphi)\subset\Spin(2,n)$ fixes a spinor both in
$\Delta_V$ and
$\Delta_W$ and this proves that we can find parallel
spinors on the cones $\hat{H}$ and $\hat{L}$.
The parallel spinor  
spinor on $\hat{H}$ gives rise to the Killing spinor $\psi_{H-}$
inducing
the Lorentzian Einstein Sasaki structure on $H$, as we discussed 
it already.
And now we can also see that there exists a 
Killing spinor $\psi_{L-}$ on 
the Riemannian Einstein spin space $L$ of positive scalar 
curvature. The associated vector to the spinor $\psi_{L-}$
on the Riemannian space $L$ vanishes, since the symplectic form 
$\omega_o$ lives on $V$ only. All together,   
this leads us to the following characterisation 
result. Thereby, we say that a conformal 
Killing spinor $\varphi_-$ has no singularities on $M$ 
if it has no zero and the function $\|\alpha_{\varphi_-}\|^2=
g(\alpha_{\varphi_-},\alpha_{\varphi_-})$ either has no zero
or is identically zero. 

\begin{THEO} \label{TH5} Let $\varphi_-$ be a conformal Killing spinor
without singularities on a simply connected Lorentzian spin manifold 
$(M^n,g)$ of dimension
$n$. Let $\alpha_{\varphi-}$ be the dual of the Dirac current
with rank $rk(\alpha_{\varphi-})$ and length 
$\|\alpha_{\varphi-}\|^2$. The following cases occur.  
\begin{enumerate}
\item
It is \[rk(\alpha_{\varphi-})=0\qquad and 
\qquad\|\alpha_{\varphi-}\|^2=0\]
and then $\varphi_-$ is locally conformally equivalent to a
parallel spinor on a Brinkmann space with lightlike Dirac current.
\item It is \[rk(\alpha_{\varphi-})=0\qquad and\qquad 
\|\alpha_{\varphi-}\|^2<0\]
and then $[g]=[-dt^2+h]$, where $h$ is a Ricci-flat Riemannian metric
admitting a parallel spinor.
\item The dimension $n$ is odd and the rank 
\[rk(\alpha_{\varphi-})=
(n-1)/2\] is maximal. Then $\varphi_-$ is conformally equivalent to
a Killing spinor on a Lorentzian Einstein-Sasaki manifold.
In this scaling the Dirac current is timelike of constant length.
\item The dimension $n$ is even and \[rk(\alpha_{\varphi-})=(n-2)/2\ .\]
Then it is $(d\alpha_{\varphi-})^{n/2}\neq 0$ and 
$(M,g)$ is a Fefferman
space with twistor spinor $\varphi_-$ and $\|\alpha_{\varphi-}\|^2=0$.
\item It is \[0<rk(\alpha_{\varphi-})<(n-2)/2\] and then there 
is a product metric $h\times l$ in $[g]$, whereby $h$ is 
an Einstein-Sasaki metric on a Lorentzian
space $H$ of dimension $p:=2\cdot rk(\alpha_{\varphi-})+1$ admitting
a Killing spinor $\psi_{H-}$ and $l$ is an Einstein metric with Killing 
spinor $\psi_{L-}$, whose associated $1$-form is trivial, 
on the Riemannian space $L$ with positive scalar curvature
$scal_l=-\frac{(n-p)(n-p-1)}{p(p-1)}scal_h$. The spinor 
$\psi_{H-}\otimes\psi_{L-}\in\Gamma(S)$ is a twistor spinor on $M$ with
timelike Dirac current (of constant negative length in the scaling $h\times 
l$).
\end{enumerate}
\end{THEO}
We want to make a remark to the last point of Theorem \ref{TH5}. 
Indeed,
the formula (\ref{For1})  for the holonomy shows in general
that if 
$g=h\times l$
is a product metric such that $h$ and $l$ admit Killing spinors
$\psi_{H-}$ and $\psi_{L-}$ with Killing numbers 
$\lambda_H=\pm i\lambda_L$ then the tensor product
$\psi_{H-}\otimes\psi_{L-}\in\Gamma(S^g)$ is a conformal Killing spinor
on $h\times l$.  One can also prove this fact by confirming the twistor 
equation
directly using the spinor connection
\[\nabla^{S,g}=\nabla^{S,h}\otimes\1\ +\ \1\otimes\nabla^{S,l}\ .\]
However, for this one has to work out carefully the correct 
identification
for an appropriate tensor product of the Clifford algebras 
$Cl(1,p-1)$ or $Cl(p-1,1)$ with $Cl(0,n-p+1)$ or $Cl(n-p+1,0)$
on the one side 
and the Clifford algebra
$Cl(1,n-1)$ or $Cl(n-1,1)$ on the other side.

%%%%%%%%%%%%%%%%%%%%%%%%%%%%%%%%%%%%%%%%%%%%%%%%%%%%%%%%%%%%%
%\section{Solutions with singualarities}%%% 
%%%%%%%%%%%%%%%%%%%%%%%%%%%%%%%%%%%%%%%%%%%%%%%%%%%%%%%%%%%%%

%%%%%%%%%%%%%%%%%%%%%%%%%%%%%%%%%%%%%%%%%%%%%%%%%%%%%%%%%%%%%
%\section{Normal conformal holonopy description in dimension 4}%%% 
%%%%%%%%%%%%%%%%%%%%%%%%%%%%%%%%%%%%%%%%%%%%%%%%%%%%%%%%%%%%%

%%%%%%%%%%%%%%%%%%%%%%%%%%%%%%%%%%%%%%%%%%%%%%%%%%%%%%%%%%%%%
\section{Further questions and outlook  }%%%
%%%%%%%%%%%%%%%%%%%%%%%%%%%%%%%%%%%%%%%%%%%%%%%%%%%%%%%%%%%%%
\label{p11}

We are concerned in this paper with the study of solutions of
the normal twistor equations in conformal geometry
and their relation to the normal conformal holonomy 
representation.
The discussion shows that for decomposable and
weakly irreducible holonomy representations without 
dilation
the conformal geometry of the underlying space of a solution can be 
described by (products of) special geometries on Einstein spaces
or at least Ricci-isotropic spaces.
Those geometric structures are well-known and were subject to substantial
work in the literature in the past
and, therefore, there is
considerable knowledge that could be applied here for 
further and more particular interest in those solutions. Moreover, there 
is
a well-known case in the literature, which has in
our context an irreducible normal conformal
holonomy representation. We mean the Fefferman spaces
in pseudo-Riemannian geometry. However, at this point 
we already come across a quite natural question 
that does not seem to have a complete answer.\\

{\it
What is the geometric structure of a conformal space
that has irreducible normal conformal holonomy representation?}\\

Or more general:\\

{\it How does the complete list of possible normal conformal holonomy 
groups look like and what kind of conformal structures do they 
describe?}\\

We want to mention two particular examples to illustrate the 
questions that we ask for.
First, there is the exceptional group
$G^*_2\subset\SO(3,4)$ in the split case of dimension $7$.
The group $\SO(3,4)$ is the M{\"o}bius group  for
a conformal space in signature $(2,3)$. The first question then 
asks for the underlying geometry (and their existence) 
that belongs to 
a conformal space, which has normal conformal holonomy group
$G_2^*$. 

A further example of an interesting holonomy
representation according to the second
(more general) question is the case of a weakly
irreducible holonomy group with dilation in 
the M{\"o}bius group $\SO(1,n+1)$, which then belongs
to conformal Riemannian geometry. (Those spaces do not
have solutions for the normal twistor equations.)
One can ask this question 
in the first instance for $3$-dimensional Riemannian space. It seems
that a normal conformal holonomy classification  
even in this case is not completely known.
 
In general, conformal Killing forms were already
introduced and studied in the works \cite{Kas68}, \cite{Tac69}.  
Recently,
there was a systematic investigation on this topic by
U. Semmelmann (cf. \cite{Sem01}). In particular, this
work shows the construction of a
certain Killing connection,
with respect to whom all conformal Killing forms
find an interpretation as parallel sections
in a certain differential form bundle.
(In fact, it is the same bundle as our M{\"o}bius form
bundle.) It arises now the question what properties 
does this Killing connection have, in particular, which 
structure group is attached to it, and whether this connection
is somehow related to a conformally invariant connection
(e.g. the normal one) on the M{\"o}bius frame bundle
with structure group $\SO(r+1,s+1)$. 

An extension of our investigations that leads 
in a similar direction, as proposed for the Killing connection
above,
is the idea 
of dropping the 
normalisation 
condition
for our twistor equations and do investigations for more general
conformal
Killing forms. For such
an attempt one could ask what conformally invariant
connections on $\frak{M}(M)$, that induce the twistor
equations, are reasonable to 
look at? For example,
does there always exists 
a conformal connection with
structure group in $\SO(r+1,s+1)$, which induces appropriate twistor
equations to a given conformal Killing form, or what
other kind of structure groups should be considered.
However, these ideas are bit loose here for the moment. Nevertheless,
there is the question whether the existence of
a conformal Killing form gives rise to some twistor
with a certain stabiliser group and some other tensor(s), which
describe uniquely the underlying conformal geometry of a space,
and are useful in order to find a systematic construction
principle for those solutions.

%-----------------------------------------------------------------------

\end{sloppypar}
\end{document}